\documentclass[preprint,11pt]{elsarticle}
\usepackage{amssymb,amsmath,float,graphicx,appendix,float}
\usepackage{subfigure}
\usepackage{url,color}
\usepackage{afterpage}

\newcommand{\bsub}{\begin{subequations}}
\newcommand{\esub}{\end{subequations}$\!$}

\newcommand{\ds}[0]{\displaystyle}

\newcommand{\eps}{{\displaystyle\varepsilon}}

\newcommand{\jlb}[1]{\textcolor{black}{#1}}

\newcommand{\jl}[1]{\textcolor{black}{#1}}
\newcommand{\jla}[1]{\textcolor{black}{#1}}
\newcommand{\al}[1]{\textcolor{black}{#1}}
\newcommand{\ala}[1]{\textcolor{black}{#1}}

\newcommand{\kg}[1]{\textcolor{black}{#1}}

\newtheorem{thm}{Theorem}[section]

\numberwithin{coronumber}{thm}

\newtheorem{proof}{Proof}[section]

\numberwithin{equation}{section}

\newcommand{\bigoh}{\mathcal{O}}
\newcommand{\littleoh}{ \mbox{{\scriptsize $\mathcal{O}$}}}

\begin{document}

\begin{frontmatter}

\journal{}

\title{Derivation and Equilibrium Analysis of a Regularized Model for Electrostatic MEMS.}

\author{A. E. Lindsay}
\ead{a.lindsay@nd.edu}
\address{Department of Applied and Computational Mathematics and Statistics,\\
University of Notre Dame, South Bend, Indiana, 46556, USA.}

\author{J. Lega}
\ead{lega@math.arizona.edu}
\address{Department of Mathematics, University of Arizona, Tucson, Arizona, 85721, USA.}

\author{K. B. Glasner}
\ead{kglasner@math.arizona.edu}
\address{Department of Mathematics, University of Arizona, Tucson, Arizona, 85721, USA.}

\begin{abstract}
In canonical models of Micro-Electro Mechanical Systems (MEMS), an event called touchdown whereby the electrical components of the device come into contact, is characterized by a blow up in the governing equations and a non-physical divergence of the electric field. In the present work, we derive novel regularized governing equations whose solutions remain finite at touchdown and exhibit additional dynamics beyond this initial event before eventually relaxing to new stable equilibria. \ala{We employ techniques from variational calculus, dynamical systems and singular perturbation theory to obtain a detailed understanding of the novel behaviors exhibited by the regularized family of equations.}
\end{abstract}


\end{frontmatter}

\section{Introduction and statement of main results.}

Micro-Electro Mechanical Systems (MEMS) are a large collection of miniaturized integrated circuits and moving mechanical components that can be fabricated together to perform a multitude of tasks. MEMS practitioners aim to manipulate the interaction between electrostatic forces and elastic surfaces to design a variety of complex devices with applications in every area of science and industry. In such interactions, the elastic surfaces of a MEMS device may be overwhelmed if the electrostatic forces acting on them are too strong. Such a failure in a MEMS device is manifested by an instability, known as the \emph{pull-in instability}.

\begin{figure}[htbp]
\centering
\includegraphics[height=4cm,clip]{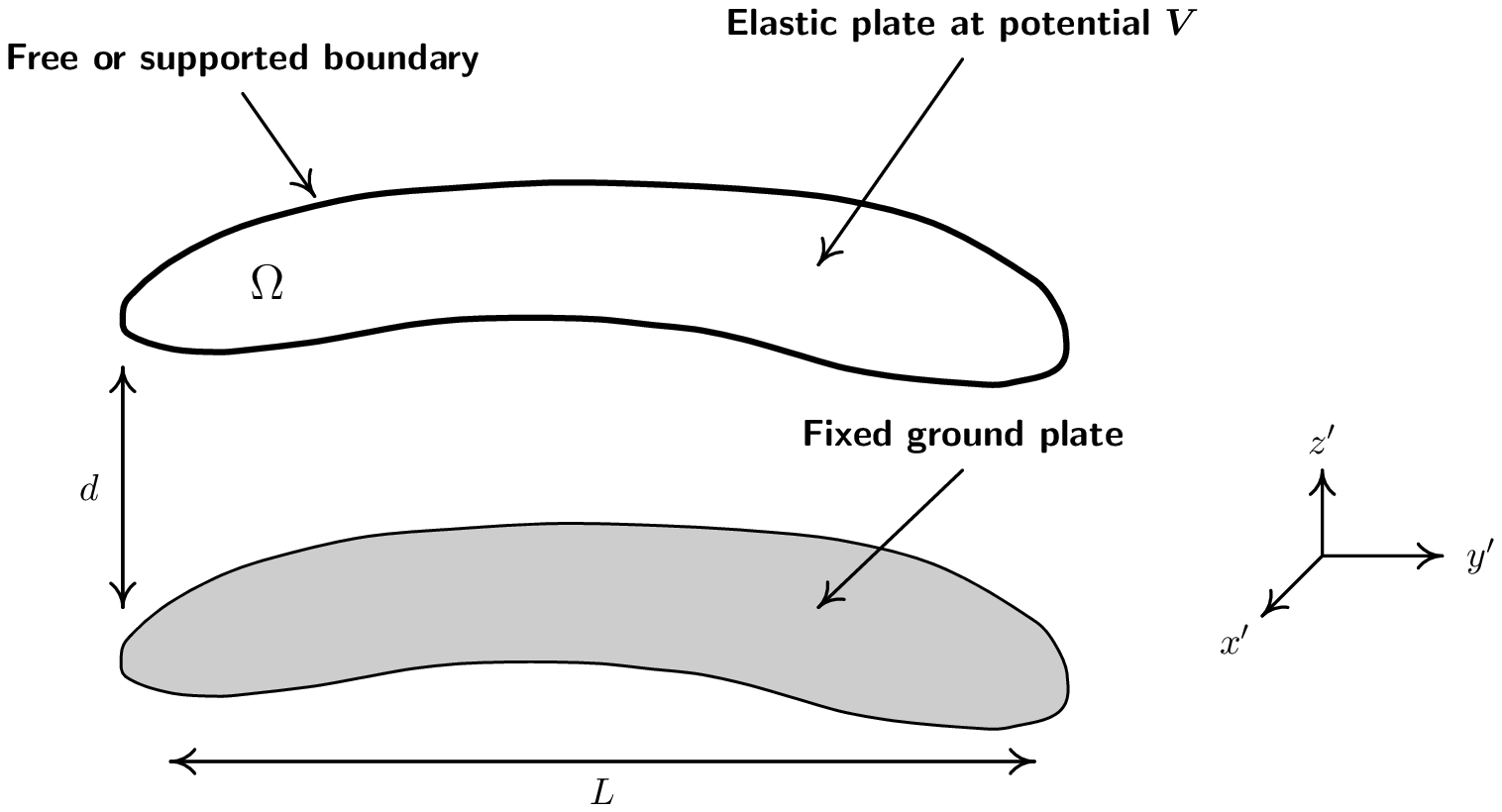}
\parbox{4.5in}{\caption{Schematic diagram of a MEMS capacitor (reproduced from \cite{LW}).\label{fig:schem}}}
\end{figure}

In a capacitor type MEMS device, an elastic membrane is held fixed along its boundary above an inelastic substrate. When an electric potential $V$ is applied between these surfaces, the upper elastic surface deflects downwards towards the substrate. If $V$ is small enough, the deflection will reach an equilibrium, however, if $V$ exceeds the \emph{pull-in voltage} $V^{\ast}$, no equilibrium configuration is attainable and the top plate will \emph{touch down} on the substrate.  Figure \ref{fig:schem} contains a schematic representation of the device.

Touchdown is a very rapid event whereby large quantities of energy are focused on small spatial regions of the MEMS device over short time scales. Consequently this process develops large forces at specific areas which can be either useful to the operation of the device or destructive. In many mathematical models of MEMS, touchdown is captured by \jl{finite time quenching}, \emph{e.g.} blow-up of solution derivative and energy. \al{Accordingly, many important operational aspects of MEMS, such as the time and location of touchdown, can be investigated by studying this quenching event.}

However, a loss of existence to model solutions results in no information regarding configurations of MEMS after a \jl{primary} touchdown event. This paper presents an initial attempt to describe behavior of MEMS after \jl{touchdown}. To this end, we derive the second order equation
\bsub\label{mems_intro_eqns}
\begin{equation}\label{mems_intro_eqns_a}
u_t = \Delta u - \frac{\lambda}{(1+u)^2} + \frac{\lambda\eps^{m-2}}{(1+u)^m},  \quad  x\in\Omega; \qquad u = 0, \quad x\in\partial\Omega,
\end{equation}
which models the dimensionless deflection $u(x,t)$ as \jl{that of} a membrane, and the fourth order problem
\begin{equation}\label{mems_intro_eqns_b}
u_t = -\Delta^2 u - \frac{\lambda}{(1+u)^2} + \frac{\lambda\eps^{m-2}}{(1+u)^m},  \quad  x\in\Omega; \qquad u = \partial_n u = 0, \quad x\in\partial\Omega,
\end{equation}
\esub
which is a beam description of the deflecting surface. \al{The modelling literature on MEMS has involved second (cf. \cite{GUO3,GUO4,GPW}) and fourth order (cf. \cite{LLS2013,LL,LY,GUO,LaurencotWalker2}) descriptions of the elastic nature of the deflecting surface and so we aim to investigate the effects of regularization on both.} In both cases, $\Omega$ is a bounded region of $\mathbb{R}^n$ and $\lambda \propto V^{2}$ is a parameter quantifying the relative importance of elastic to electrostatic forces. The physically relevant dimensions are $n=1,2$. The small parameter $\eps$ in \eqref{mems_intro_eqns} mimics the effect of a small insulating layer placed on top of the substrate to prevent a short circuit of the device as the gap spacing $1+u$, \jlb{$u<0$,} \jl{locally} shrinks to zero.

For the case $\eps=0$, equations \eqref{mems_intro_eqns} reduce to canonical models originally introduced by Pelesko (cf. \cite{PB}), the salient properties of which are now well known. Of particular importance amongst the many results, is the existence of a pull-in voltage $\lambda^{\ast}$ such that if $\lambda<\lambda^{\ast}$, then $u(x,t)$ approaches a unique and stable equilibrium as $t\to\infty$, while for $\lambda>\lambda^{\ast}$ no equilibrium solutions are possible and $u(x,t)$ reaches $-1$ in some finite time, $t_c$. In the 1D setting, the equilibrium structure consists of one stable and one unstable branch that meet at $\lambda^{\ast}$ (cf. dashed curve of Fig.~\ref{fig:intro_bif_diagrams}). In the case where $\lambda>\lambda^{\ast}$, there have been many studies centred on describing the local properties of the device near touchdown. For example, in the second order equation,
\bsub\label{memb_intro}
\begin{align}
\label{memb_intro_a} u_t = \Delta u - \frac{\lambda}{(1+u)^2}, \quad x\in  \Omega,\\[5pt]
\label{memb_intro_b} u(x,0) = u_0(x), \quad x\in\Omega; \qquad u = 0, \quad x\in\partial\Omega,
\end{align}
\esub
a detailed analysis \cite{GPW} of solutions near touchdown revealed the local behavior
\begin{equation}\label{mems_intro_2}
 u \to -1 + [3\lambda(t_c-t)]^{1/3}\Big( 1 -\frac{1}{2|\log(t_c-t)|} + \frac{(x-x_c)^2}{4(t_c-t)|\log(t_c-t)|} +\cdots\Big),
\end{equation}
in the vicinity of the touchdown point $x_c$, for $t\to t_c^{-}$. Detailed scaling laws for $t_c$ in the limits $\lambda\to\infty$ and $\lambda-\lambda^{*}\to0^{+}$ have also been established in \cite{GUO3, GUO4}. In the fourth order problem,
\bsub\label{biharm_intro}
\begin{align}
\label{biharm_intro_a} u_t = -\Delta^2 u - \frac{\lambda}{(1+u)^2}, \quad x\in  \Omega,\\[5pt]
\label{biharm_intro_b} u(x,0) = u_0(x), \quad x\in\Omega; \qquad u = 0, \quad \partial_n u = 0, \qquad x\in\partial\Omega,
\end{align}
\esub
less is known about the equilibrium structures and dynamics of touchdown in the absence of static solutions. In the special cases where $\Omega$ is the unit strip $[-1,1]$ or the unit disc $\{ x\in\mathbb{R}^2 \, \mid \, |x| \leq1 \}$, the existence of the pull-in voltage $\lambda^{\ast}$ was shown in \cite{LW}. Similar results were obtained in \cite{GUO} for the case where pinned boundary conditions $u=\Delta u =0$ were applied to the boundary. For $\lambda>\lambda^{\ast}$ and for $\Omega$ the unit strip $[-1,1]$ or the unit disc $\{x\in\mathbb{R}^2 \, \mid\, |x|\leq1\}$, it was shown in \cite{LL} that the device touches down in finite time $t_c$. A detailed numerical and asymptotic study established the local behavior
\begin{equation}\label{intro:biharm_TD}
u(x,t) \to -1 + (t_c-t)^{1/3} v(y), \qquad  y = \frac{x-x_c}{(t_c-t)^{\frac14}}\jla{\lambda^{1/4}}, \qquad t \to t_c^{-},
\end{equation}
where $v(y)$ is a self-similar profile satisfying an associated ordinary differential equation. In addition to the local behavior of solutions as $t\to t_c^{-}$, the fourth order problems \eqref{biharm_intro} have additional interesting dynamical features whereby touchdown can occur simultaneously at multiple points of the domain. In one dimension \cite{LL}, the singularities can form at two distinct points separated about the origin. In two dimensions \cite{LLS2013}, the multiplicity of singularities can be greater with the exact quenching set depending delicately on the geometry of the boundary and the parameter $\lambda$.

The rich dynamical behavior associated with the touchdown event raises the interesting question of how one can make sense of solutions to \eqref{memb_intro} and \eqref{biharm_intro}, and understand \jl{the} behavior of MEMS after touchdown. The finite time singularities exhibited by \eqref{memb_intro} and \eqref{biharm_intro} result in the gap spacing $1+u$ becoming arbitrarily small as $t\to t_c^{-}$ for $\lambda$ sufficiently large. Consequently, a physically unreasonable situation occurs -  the electric field generated between the plates becomes arbitrarily large as $t\to t_c^{-}$. The focus of this paper is first to regularize the singularity in the electric field at touchdown, thereby rendering it large but finite thereafter, and second to describe the post-touchdown equilibrium configurations of the resulting model. We derive suitable regularized equations in Section \ref{sec:derivation} and \jlb{analyze their properties in Section \ref{sec:Prop}. First, we} show in Section \ref{sec:WP} that the regularized equations are globally well-posed. The variational nature of these equations \jlb{then} leads us to consider equilibrium solutions. Numerical simulations \jlb{shown in Section \ref{sec:var_dyn}} indicate that the regularized equations we propose undergo additional dynamics beyond the initial touchdown event (see for instance Fig.~\ref{fig:intro_spreading}) and converge towards \jlb{a new branch of} equilibrium solutions. \jlb{We show the corresponding bifurcation diagrams in Section \ref{sec:bif_diagrams} and explain how the new branch of solutions appears in Section \ref{sec:new_sol}.} We then describe the properties of post-touchdown equilibrium configurations in terms of matched asymptotic expansions \jlb{in Section \ref{sec:scalings}}. We summarize our results in Section \ref{sec:discussion} and discuss implications of the present work, in particular regarding the bistable nature of the proposed regularized equations.

\section{Derivation of regularized governing equations.}\label{sec:derivation}

In this section we derive a new model for the operation of a MEMS device \jl{with} a small insulating layer resting on the substrate\jl{,} whose purpose \jl{is} to physically prevent \jl{the occurrence of} a short circuit. Based on this principle, the new model features an obstacle type regularization of touchdown\jl{, in the form of} a perturbed electrostatic potential with a repulsive term \jl{that} mimics the obstacle.

In dimensional form, the model requires that the \jl{vertical (ie. parallel to the $z$-direction)} deflection $u(x,y,t)$ of a plate occupying a region $\Omega\subset\mathbb{R}^2$ with boundary $\partial\Omega$, satisfies \cite{PB}
\bsub\label{memsdim}
\begin{align}
\label{memsdima}
\rho h\, \frac{\partial^2 u}{\partial t^2} + a\,\frac{\partial u}{\partial t} + EI \, \Delta^2_{\perp} u - T\, \Delta_{\perp} u &= -\frac{\epsilon_0}{2}\,|\nabla\phi|^2_{z=u}  \qquad x \in\Omega;\\[5pt]
\label{memsdimb} \nabla\cdot(\sigma \nabla \phi) &= 0 \qquad -(d +h) \leq z  \leq u(x,y,t),
\end{align}
where $\perp$ indicates differentiation with respect to the $x$ \jl{and} $y$ directions, and the \jl{permittivity} $\sigma$ satisfies
\begin{equation}\label{memsdimc}
\sigma = \left\{  \begin{array}{rl} \sigma_0, & -d \leq z \leq u(x,y,t)\\[5pt] \sigma_1, & -(d+h) \leq z \leq -d \end{array} \right. .
\end{equation}
\al{In equations \eqref{memsdim}, $\rho$ $h$, $EI$ and $T$ are the density per unit length, thickness, flexural rigidity and tensile load of the plate. The parameter $a$ represents damping forces on the system, $\epsilon_0$ is the permittivity of free space and $d$ is the undeflected gap spacing.} The \jl{electric} potential \jl{$\phi$} at the ground plate is zero and a voltage $V$ is applied on the upper plate so that
\begin{equation}\label{memsdimd}
\phi(-(d+h)) = 0,\qquad \phi(u) = V.
\end{equation}
\esub
The problem is now reduced by recasting equations \eqref{memsdim} in the dimensionless variables
\[
x' = \frac{x}{L}\quad y' = \frac{y}{L} \quad z' = \frac{z}{d}, \quad u' = \frac{u}{d},\quad \phi'=\frac{\phi}{V}, \quad \sigma' =\frac{\sigma}{\sigma_0}
\]
and applying the small aspect ratio $\delta\equiv d/L \ll1$. \jl{Here, $L$ is a characteristic linear dimension of the domain $\Omega$.} Concentrating first on the potential equation \eqref{memsdimb}, the non-dimensional equation \jl{for} $\phi'$ satisfies
\bsub\label{nondim1}
\begin{align}
\label{nondim1a} \nabla'\cdot(\sigma' \nabla' \phi' ) &= 0, \qquad -(1 +h/d) \leq z'  \leq u'(x',y',t);\\[5pt]
\label{nondim1b} \sigma' &= \left\{  \begin{array}{cl} 1, & -1 \leq z' \leq u'(x',y',t);\\[5pt] \ds\frac{\sigma_1}{\sigma_0}, & -(1+h/d) \leq z' \leq -1 \end{array} \right.\\[5pt]
\phi'(-(1+h/d)) &= 0,\qquad \phi'(u') = 1.
\end{align}
\esub
In non-dimensional coordinates, we have that
 \[ \nabla' \equiv \left( \frac{1}{L} \frac{\partial}{\partial x'}, \ \frac1L\frac{\partial}{\partial y'} ,\ \frac{1}{d} \frac{\partial}{\partial z'}\right)\]
and therefore problem \eqref{nondim1} reduces to
\bsub\label{nondim2}
\begin{align}
\label{nondim2a}\frac{\partial^2 \phi^{'}_{+}}{\partial z'^2} + \delta^2 \left(\frac{\partial^2 \phi^{'}_{+}}{\partial x'^2} +\frac{\partial^2 \phi^{'}_{+}}{\partial y'^2} \right) &=0, \qquad -1\leq z' \leq u'; \\[5pt]
\label{nondim2b}\frac{\partial^2 \phi^{'}_{-}}{\partial z'^2} + \delta^2 \left(\frac{\partial^2 \phi^{'}_{-}}{\partial x'^2} +\frac{\partial^2 \phi^{'}_{-}}{\partial y'^2} \right) &=0, \qquad -1-\frac{h}{d}\leq z' \leq -1;\\[5pt]
\label{nondim2c} \phi^{'}_{+}(u') = 1, \qquad \phi^{'}_{+}(-1) =\phi^{'}_{-}(-1), & \qquad \frac{\partial}{\partial z'}\phi^{'}_{+}(-1) = \frac{\sigma_1}{\sigma_0}\frac{\partial}{\partial z'}\phi^{'}_{-}(-1), \qquad \phi^{'}_{-} (-1-d/h) = 0.
\end{align}
\esub
Applying the small aspect ratio $\delta\to0$, the leading order solution to \eqref{nondim2} is
\begin{equation}\label{nondim3}
\phi' = \left\{ \begin{array}{cc} 1 + \ds\frac{z'-u'}{(1+u') + \ds\frac{d\sigma_0}{h\sigma_1}} & \jla{-}1\leq z' \leq u'; \\[12pt]
\ds\frac{z'+1+\ds\frac{d}{h}}{\ds\frac{\sigma_1}{\sigma_0}(1+u') + \ds\frac{d}{h}}&-1-\ds\frac{h}{d} \leq z' \leq -1 \end{array}\right.
\end{equation}
The explicit solution \eqref{nondim3} which arises in this small aspect ratio limit affords a significant reduction in the complexity of the governing equations. If the limit $\delta\to0$ is not exercised, the system for the potential \eqref{nondim2} and the non-dimensionalized form of \eqref{memsdima} constitute a free boundary problem for the deflection $u(x,y,t)$ of the device. With the exclusion of the insulating layer introduced here in \eqref{memsdimc}, the qualitative properties of dynamic and steady solutions of this free boundary problem have been studied in \cite{LaurencotWalker1,LaurencotEscherWalker1,LaurencotEscherWalker2}. These studies have established the well-posedness theory for the system of evolution equations \eqref{memsdim}, the existence of a pull in voltage and also the convergence of equilibrium solutions of the free boundary problem to those of the small aspect limit as $\delta\to0$. Accordingly, there is good reason to believe that the small aspect ratio approximation is a good one. In light of the significant simplifications it affords, we proceed by calculating from \eqref{nondim3} that the forcing on the surface $z'=u'(x',y')$ is given by
\begin{equation}\label{nondim4}
\frac{\epsilon_0}{2}|\nabla \phi|^2 = V^2\frac{\epsilon_0}{2 \jla{d^2}} \left[ \left(\frac{\partial\phi'}{\partial z'}\right)^2 +\mathcal{O}(\delta^2)   \right ] = V^2\frac{\epsilon_0}{2 \jla{d^2}}\ds\frac{1}{\Big(1+u' + \ds\frac{d\sigma_0}{h\sigma_1}\Big)^2}.
\end{equation}
After selecting the time scale $t = (L^2a/T)t'$ in \eqref{memsdim} and substituting the reduced term arrived at in \eqref{nondim4}, the equation
\bsub\label{nondim5}
\begin{equation}\label{nondim5a}
\alpha^2\, \frac{\partial^2 u'}{\partial t'^2} + \,\frac{\partial u'}{\partial t'} +\,\beta \Delta_{\perp}^{'2} u'  - \, \Delta^{'}_{\perp} u' = -\frac{\lambda}{(1+u' + \eps)^2}\\[5pt]
\end{equation}
is obtained, where the dimensionless groups are
\begin{equation}\label{nondim5b}
\beta = \frac{EI}{L^2T}, \qquad \alpha = \frac{\sqrt{T\rho h}}{aL},\qquad \eps =\frac{d\sigma_0}{h\sigma_1}, \qquad \lambda=\frac{\epsilon_0 L^2V^2}{2d^3T}.
\end{equation}
\esub
The focus of our attention is further restricted to the case of small quality factor for which the $\alpha^2 u_{tt}$ term in (\ref{nondim5}) is considered negligible. This approximation, called the viscous damping limit \cite{PB}, assumes that inertial effects are negligible compared to those of damping. All quantities are now dimensionless and all derivatives are in the $x, y$ directions so the $'$ and $\perp$ notations can be dropped. \jl{In summary, the dynamics of a MEMS device in the presence of an insulating layer is thus modeled by the following} obstacle problem
\bsub\label{mems_main}
\begin{align}
\label{mems_main_a} u_t = -\beta \Delta^2 u + \Delta u - \jl{\frac{d \psi_\eps}{d u}, \qquad \psi_\eps(u) = - } \frac{\lambda}{(1+u+\eps)}, \quad  x\in \Omega;\\[5pt]
\label{mems_main_b} u \geq -1, \qquad  x\in \Omega;
\end{align}
with boundary and initial values
\begin{equation}\label{mems_main_c}
u = 0, \quad \partial_n u = 0, \quad \mbox{on}\quad \partial\Omega; \qquad u = 0, \quad t=0.
\end{equation}
\esub
The combination of the $\eps$ term in the Coulomb nonlinearity of \eqref{mems_main_a} and the obstacle constraint \eqref{mems_main_b}, act to prevent blow up at touchdown.

\subsection{\kg{Variational nature of the obstacle problem and a regularization}}\label{sec:VI}

Obstacle problems like \eqref{mems_main} often arise in mechanics when constraints are present
\cite{kikuchioden}.  These problems are commonly written as variational inequalities, whose
basic mathematical properties such as existence and uniqueness are well-established
(e.g. \cite{duvantlions1972,kinderlehrerstampacchia}).
The evolution equation \eqref{mems_main} can be thought formally as the $L^2$ gradient flow of the functional
$E:H^2(\Omega) \to \mathbb{R} \cup \{+\infty\}$ given by
\begin{equation}
  \label{obstacle}
E = \int_{\Omega} \frac{\beta}{2} (\Delta u)^2 + \frac12 |\nabla u|^2 + \psi(u) \, dx,
\end{equation}
where
\begin{equation}
  \label{obstaclepotential}
\psi(u; \epsilon) = \begin{cases}
-\ds\frac{\lambda}{1 + u + \epsilon} & u \ge 0 \\[5pt]
+ \infty &  u < 0.
\end{cases}
\end{equation}
The assignment of infinite energy to values of $u < -1$ encodes the obstacle constraint.

For practical purposes, it is often useful to work with a regularized version of the obstacle problem
which has smooth solutions (e.g. \cite{lions1969quelques,scholz1986}).  This typically involves,
in essence, replacing an energy functional like \eqref{obstacle} with one which is smooth but otherwise
mimics the penalization associated with the obstacle.

For our problem, we will replace the potential \eqref{obstaclepotential} with one which has the same qualitative
structure.  Specifically, the new potential $\phi_{\epsilon}$ will behave like $\psi$ in the following ways:
\begin{enumerate}
\item
For fixed values of $u> -1$, $\phi_{\epsilon}(u) \sim \psi(u; \epsilon)$ as $\epsilon \to 0$.
\item
$ \lim_{u \to -1^+} \, \phi_{\epsilon}(u) = + \infty$
\item
The value of $\psi$ which occurs at the obstacle value $u = -1$ is the same as the minimum of
$\phi_{\epsilon}(u)$.
\end{enumerate}

\begin{figure}[H]
\centering
\includegraphics[width=3in]{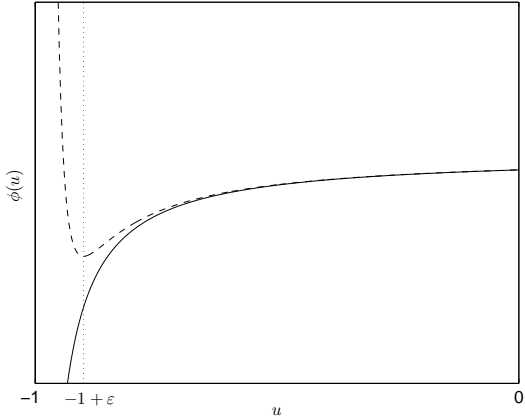}
\parbox{5.5in}{\caption{A schematic diagram of the potential \eqref{intro_potential}. The solid line indicates the case $\eps=0$ while the dashed line represent\jl{s the case} $0<\eps<1$. Note that the perturbed potential has the generic features of \jl{having} a local minimum at $u=-1+\eps$, \jl{of being} repulsive when $-1<u<-1+\eps$, and attracting when $u>-1+\eps$. \label{intro:LJ}}}
\end{figure}

A class of potentials which fulfills these criteria is
\begin{equation}
  \label{intro_potential}
\phi_{\eps}(u) = -\frac{\lambda'}{(1+u)} + \frac{\lambda' (\alpha \eps)^{m-2}}{(m-1)(1+u)^{m-1}}, \jl{\qquad \lambda' > 0,} \qquad 0<\eps<1,
\end{equation}
for integer exponents $m> 2$, and $\alpha = (2-m)/(m-1)$.  We hereafter set $\eps' = \alpha \eps$ and drop the prime.
A schematic diagram \jl{of the graph of $\phi_\eps$} is shown in Fig.~\ref{intro:LJ}.

In the preceding derivation, we have used \jl{an} elastic model of the deflecting surface based on a plate under tension, which results in a combination of Laplacian and bi-Laplacian terms in \eqref{mems_main_a}. Our analysis and observations indicate that whenever these two terms appear, the bi-Laplacian terms dominates qualitative solution features. To effect a cleaner quantitative analysis, we therefore study the bi-Laplacian and Laplacian terms in isolation, rather that in combination.
\kg{In the bi-Laplacian case we can dispense with the parameter $\beta$ by a different non-dimensionalization}
\begin{equation}\label{lambda_biharm}
 \lambda = \frac{\eps_0 L^4 V^2}{2d^3 EI}, \qquad \jl{t = \frac{L^4 a}{E I}\ t^\prime},
\end{equation}
whereas for the Laplacian case, the scaling of $\lambda$ is as in \eqref{nondim5b}.

The culmination of the obstacle regularization and separation of the linear term leads us to
study two problems, the second order equation
\bsub\label{mems_intro_3}
\begin{equation}\label{mems_intro_3a}
u_t = \Delta u- \frac{\lambda}{(1+u)^2} + \frac{\lambda\eps^{m-2}}{(1+u)^m},  \quad  x\in\Omega; \qquad u = 0, \quad x\in\partial\Omega,
\end{equation}
and the fourth order equation
\begin{equation}\label{mems_intro_3b}
u_t = -\Delta^2 u - \frac{\lambda}{(1+u)^2} + \frac{\lambda\eps^{m-2}}{(1+u)^m},  \quad  x\in\Omega; \qquad u = \partial_n u = 0, \quad x\in\partial\Omega.
\end{equation}
\esub
In particular, the singular limit $\eps\to0$ will receive special attention.

\section{Properties of the regularized equations}\label{sec:Prop}
\subsection{Well-posedness}\label{sec:WP}

In this section we detail the existence theory for both the Laplacian and bi-Laplacian problems,
which we write as
\bsub\label{WP:eqns}
\begin{align}
  \label{lp}
u_t &= \Delta u - \phi_{\eps}'(u), \quad  x\in\Omega; \qquad u = 0, \quad x\in\partial\Omega;
\\[5pt]
  \label{blp}
u_t &= -\Delta^2 u - \phi_{\eps}'(u)\quad  x\in\Omega; \qquad u = \partial_n u = 0, \quad x\in\partial\Omega,
\end{align}
\esub
together with the initial condition $u(x,0) = u_0(x)$. The spatial domain $\Omega \subset \mathbb{R}^n $ is assumed compact with a sufficiently smooth boundary. We note that the evolution equations are $L^2$ gradient flows. In particular, if
\bsub\label{energy}
\begin{eqnarray}
  \label{lenergy}
E_L(t) = \int_{\Omega}  |\nabla u|^2  + \phi_{\eps}(u) dx,\\
  \label{benergy}
E_B(t) = \int_{\Omega}  |\Delta u|^2  + \phi_{\eps}(u) dx,
\end{eqnarray}
\esub
it is easily shown that $d E_L/ dt \le 0$ and $d E_B/ dt \le 0$. The following results are proved for a class of potentials $\phi$ which is fairly general and for which \eqref{intro_potential} is a subset. For both equations we suppose
\begin{equation}\label{assumption}
\phi_{\eps}(u) \in C^1, \quad \text{$\phi_{\eps}(u) \ge \phi_{min}$ for $u \in (-1, \infty)$}, \quad
 \text{$\phi_{\eps}(u) < \phi_{max}$ for $u \in (-1+\eps,\infty)$ }.
\end{equation}
Additional restrictions for each equation are
\bsub\label{potential}
\begin{eqnarray}
  \label{potential1}
\phi_{\eps}'(u) <0 \quad \text{if $u \in (-1,-1+\eps)$, for equation (\ref{lp}), } \\[5pt]
  \label{potential2}
\phi_{\eps}(u) \sim c(\eps) (1 + u)^{-m+1} \quad u \to -1, \quad \text{ for equation (\ref{blp})}.
\end{eqnarray}
\esub
for constant $c(\eps)$.
\begin{thm}[Global Existence - Laplacian Case]
Suppose that the initial condition satisfies $u_0 \in C^0(\Omega)$ and $u_0 > -1$.
Then the solution for (\ref{lp}) exists for all $t>0$ and $u(x,t) > \min( \inf u_0, -1+\eps) $.
\end{thm}
\underline{\emph{Proof:}}
Let $u_{\pm}(t)$ solve the initial value problems
\begin{equation}
   \label{ivp}
\frac{ d  u_{\pm}}{d t} = -\phi_{\eps}'(  u_{\pm}), \quad u_-(0) = \inf u_0, \quad u_+(0) = \sup u_0.
\end{equation}
Conditions (\ref{assumption},\ref{potential1}) ensure that $u_{\pm}$ will exist for all $t>0$
and $u_{\pm} > -1$.  Furthermore, $u_- > \min( \inf u_0, -1+\eps)$.
Standard comparison methods for parabolic equations yield the {\it a priori} bounds
$u_-(t) \le u(x,t) \le u_+(t)$.  This guarantees that the solution will exist globally.
\endproof

\begin{thm}[Global Existence - bi-Laplacian Case]
Suppose that the initial condition satisfies $u_0 \in H^2(\Omega) \cap C^0(\Omega)$ and $u_0 > -1$.
Then the solution $u(x,t)$ of \eqref{blp} exists for all $t>0$, provided $m \ge 3$ in dimension $n=1$ and
$m > 3$ in dimension $n=2$.
\end{thm}

\underline{\emph{Proof:}}
Following \cite{bertozzigrunwiteslski}, it suffices to derive {\it a priori} pointwise bounds
on the solution.  This guarantees that the equation is uniformly parabolic and existence follows from
standard arguments. The gradient flow structure and $dE_{B}/dt\leq0$ implies that $E_B(T)-E_B(0)\leq 0$ for any $T >0$, and so
\begin{equation}
   \label{ineq}
\int_{\Omega} (\Delta u(T) )^2 dx \le \int_{\Omega}  (\Delta u_0 )^2  dx
+ \int_{\Omega}  \phi_{\eps}(u_0) dx - \int_{\Omega}  \phi_{\eps}(u(T)) dx.
\end{equation}
Since $\phi(\cdot)$ has a lower bound, it follows that $u \in H^2(\Omega)$ {\it a priori}.
The Sobolev imbedding theorem then gives $u \in C^1(\Omega)$ in dimension $n=1$ and $u \in C^{0,\alpha}(\Omega)$
in dimension $n=2$ where $0 < \alpha < 1$.  In particular there are constants $K_1$ and $K_2$, depending
only on the initial condition, so that
\begin{eqnarray}
  \label{sob1}
\|u\|_{C^1} < K_1, \quad n=1;\\[5pt]
  \label{sob2}
\|u\|_{C^{0,\alpha}} < K_2(\alpha), \quad n=2.
\end{eqnarray}

Now let $u_{\min} = \min u(T)$ be the minimum attained at a point $x_0$.
Note that inequality \eqref{ineq} implies an upper bound for $\int_{\Omega} \phi_{\eps}(u(T)) dx$.
In dimension $n=1$ it follows that there exist generic constants $K$ so that
\begin{equation}
C > \int_{\Omega} \phi_{\eps}(u(T)) dx \ge K(\eps) \int_{\Omega} (u_{min} + 1 + K_1 |x-x_0|)^{-m+1} dx \ge \mu( u_{min}+1),
\end{equation}
where
\begin{equation}
\mu( u_{min}+1) = K(\eps) \begin{cases}
-\ln ( u_{min}+1) & m = 3, \\
( u_{min}+1)^{-m+3} &  m > 3.
\end{cases}
\end{equation}
In dimension $n=2$ one similarly has
\begin{equation}
C >  K(\eps) \int_{\Omega} (u_{min} + 1 + K_2 |x-x_0|^{\alpha})^{-m+1} dx \ge \mu( u_{min}+1),
\end{equation}
where
\begin{equation}
\mu( u_{min}+1) = K(\alpha,\eps) \begin{cases}
-\ln ( u_{min}+1) & m = 1+ 2/\alpha, \\
( u_{min}+1)^{3-m} &  m > 1 + 2/\alpha.
\end{cases}
\end{equation}
In both cases, this establishes, for $\eps>0$, the lower bound $u > -1$ for all $t>0$.
\endproof

The two preceding results capture two important features of the perturbed potential system. First, for a wide range of potentials, equations \eqref{WP:eqns} mimic the effect of the obstacle constraint $u>-1$, established in \eqref{mems_main_b}. This provides confidence that the perturbed potential system qualitatively reflects the behavior of the obstacle problem \eqref{mems_main}. Second, in contrast to the $\eps=0$ case, the system is now well-posed for all $t>0$ and $\eps>0$ and no finite time singularity occurs. It is therefore relevant to investigate the limiting behaviour of equations \eqref{WP:eqns} in the limit $t\to\infty$. This long term behavior of equations \eqref{WP:eqns} is related to the minimizers of the functionals given in \eqref{energy}.

% for which we obtain the following results:
%
%\begin{thm}[Existence of Global Minimizers]
%Let $\Omega$ be a bounded domain of class $\mathbb{C}^{0,1}$. Then a global minimizer of \eqref{lenergy} exists in $H^1(\Omega)$ and a global minimizer of \eqref{benergy} exists in $H^2(\Omega)$.
%\end{thm}
%\underline{\emph{Proof:}}  In both the second order and fourth order case, the energies \eqref{energy} are bounded and so we can define for a minimizing sequence $\{ u_k \}_{k\in\mathbb{N}}$.

\subsection{\jlb{Variational dynamics}}\label{sec:var_dyn}

\begin{figure}[H]
\centering
\subfigure[Initial touchdown.]{\includegraphics[width=1.925in]{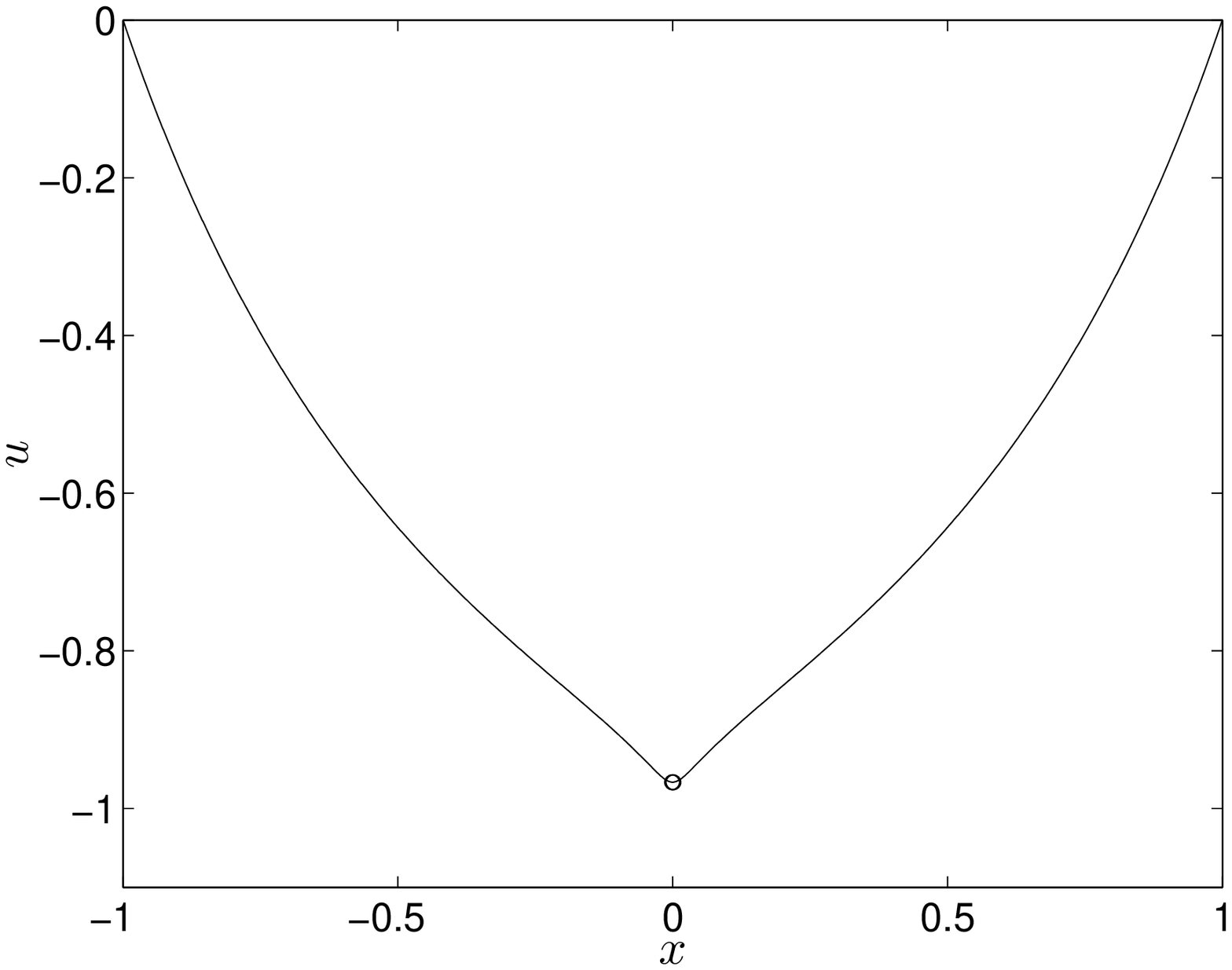}\label{fig:intro_spreading_a}}
\subfigure[Spreading of touchdown region.]{\includegraphics[width=2in]{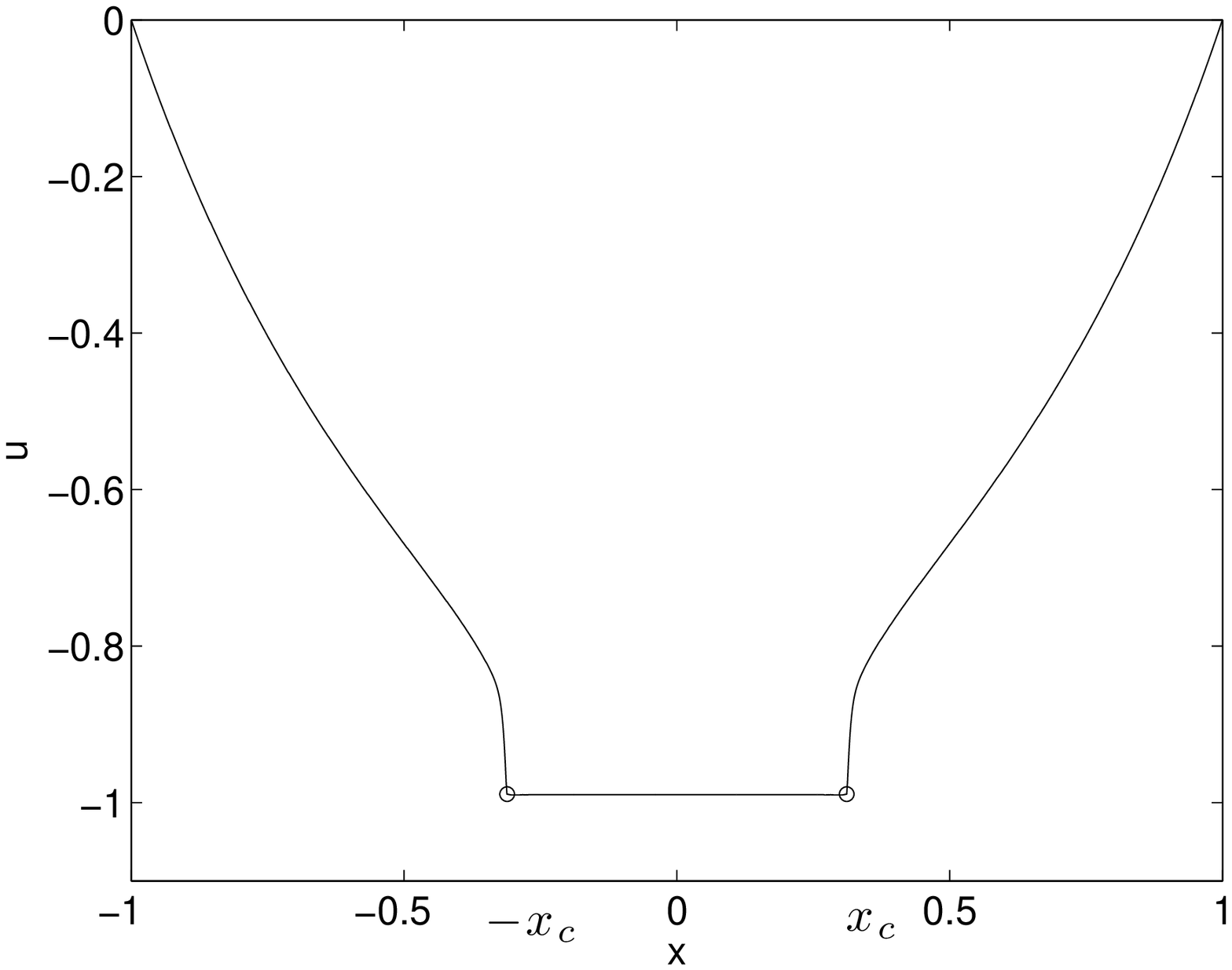}\label{fig:intro_spreading_b}}
\subfigure[Boundary pinning.]{\includegraphics[width=2in]{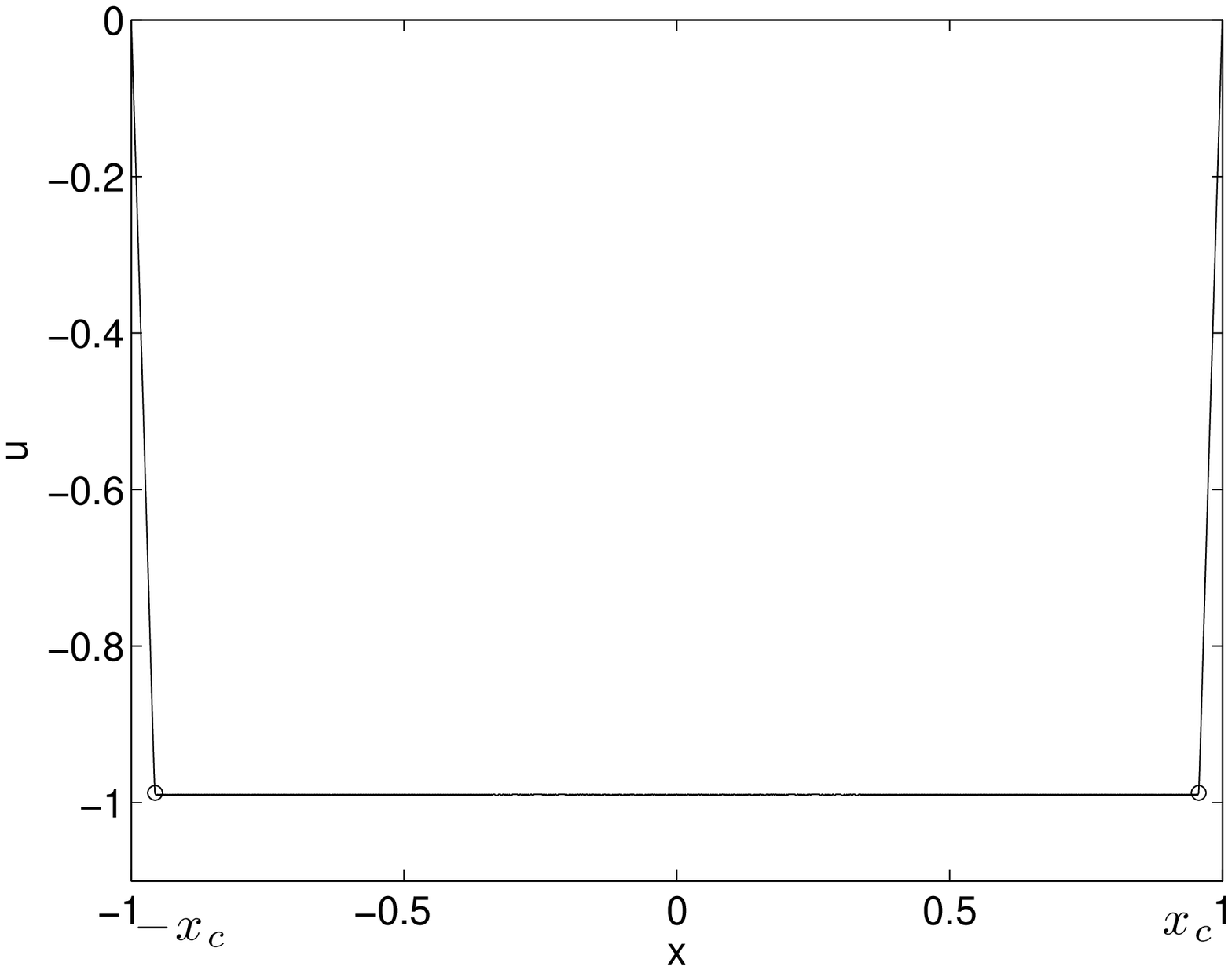}\label{fig:intro_spreading_c}}
\parbox{5.5in}{\caption{Solutions of \eqref{mems_intro_3a} initialized with zero initial data and parameter values $\eps=0.01$, $\lambda = 5$. \jl{The l}eft panel shows the initial touchdown event at $x=0$. \jl{The c}enter panel shows the spread of the touchdown region towards the boundary. Right panel: An equilibrium state is reached after the moving front is pinned by its interaction with the boundary. \label{fig:intro_spreading}}}
\end{figure}

\jl{The dynamics of Equations \eqref{mems_intro_3} is variational and leads to relaxation of the system towards equilibrium solutions. \jlb{For values of $\lambda$ such that touchdown would not occur when $\eps=0$, the regularization term in \eqref{mems_intro_3} remains of order $\eps^{m-2}$ since $1+u$ remains finite, and the dynamics in the presence of regularization is therefore a regular perturbation of the dynamics without regularization. For larger values of $\lambda$ however, the blow-up of the nonlinear term is prevented by the regularization term and the dynamics evolves towards a solution for which most of the membrane is in near contact with the dielectric layer covering the substrate. This} is illustrated in Figure \ref{fig:intro_spreading}, in the Laplacian case\jlb{, for a one-dimensional domain, $\Omega = [-1,1]$}. As an initially flat membrane deforms under the effect of the applied electric field, it first touches down at one point in the middle of the domain $\Omega$. A region where $u \simeq -1+ \eps$ then grows from the initial touchdown location towards  the boundary of the domain. The outermost points of this growing region slow down as they get closer to the edge of the domain, and eventually stop at a distance $x_c$ from the boundary. Qualitatively similar behavior is observed in the case of the bi-Laplacian.}

\jlb{This dynamics is markedly different from the $\eps=0$ case, for which no equilibrium solutions exist above a given threshold $\lambda > \lambda^{\ast}$. As we will see below, this is due to the appearance of a branch of equilibrium solutions of large $L^2$ norm, which exists when $\eps \ne 0$.}

\subsection{One-dimensional equilibrium solutions and bifurcation diagrams} \label{sec:bif_diagrams}

\begin{figure}[H]
\centering
\subfigure[Laplacian bifurcation \jl{diagram}]{\includegraphics[width=0.45\textwidth]{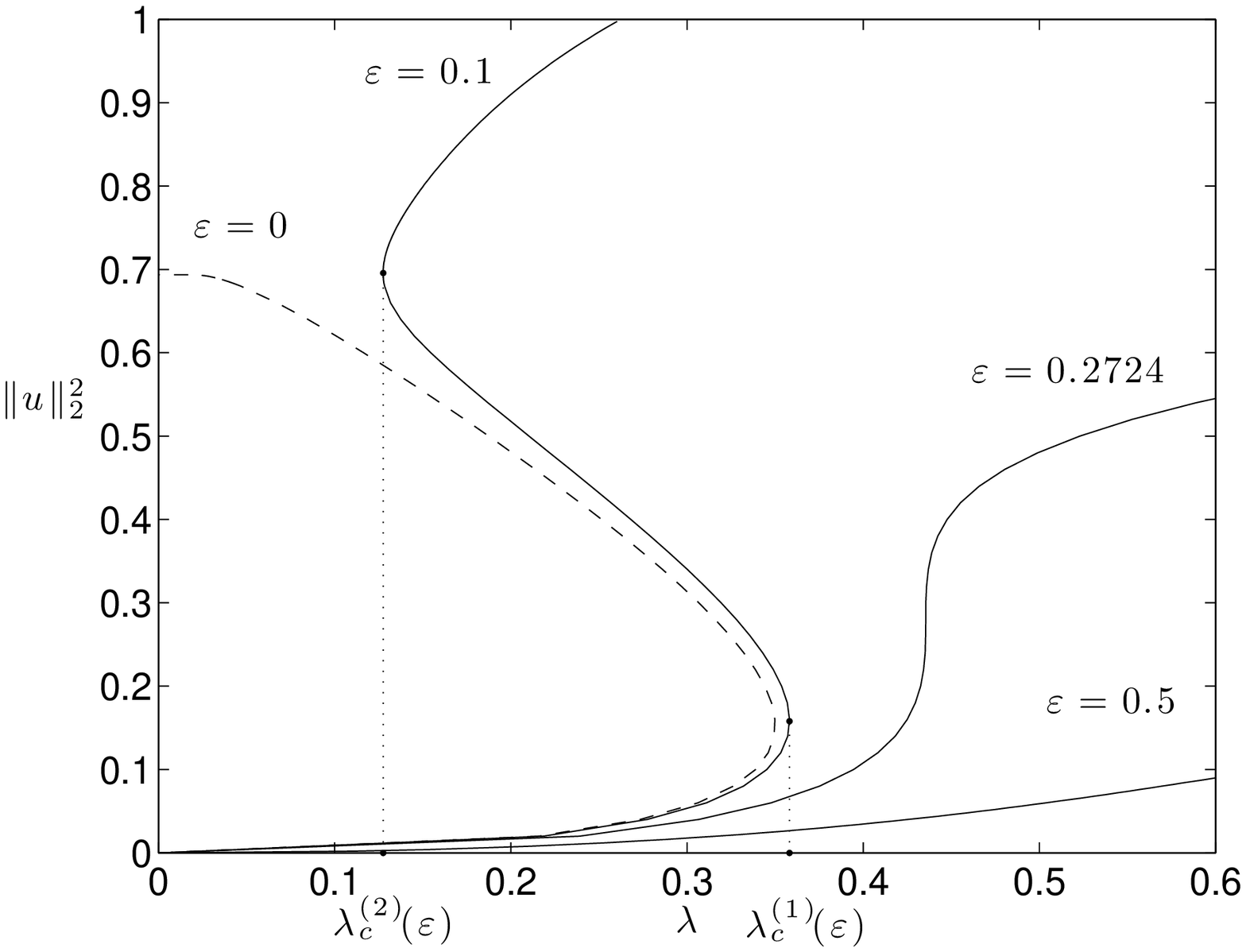}}\qquad
\subfigure[Bi-Laplacian bifurcation \jl{diagram}]{\includegraphics[width=0.45\textwidth]{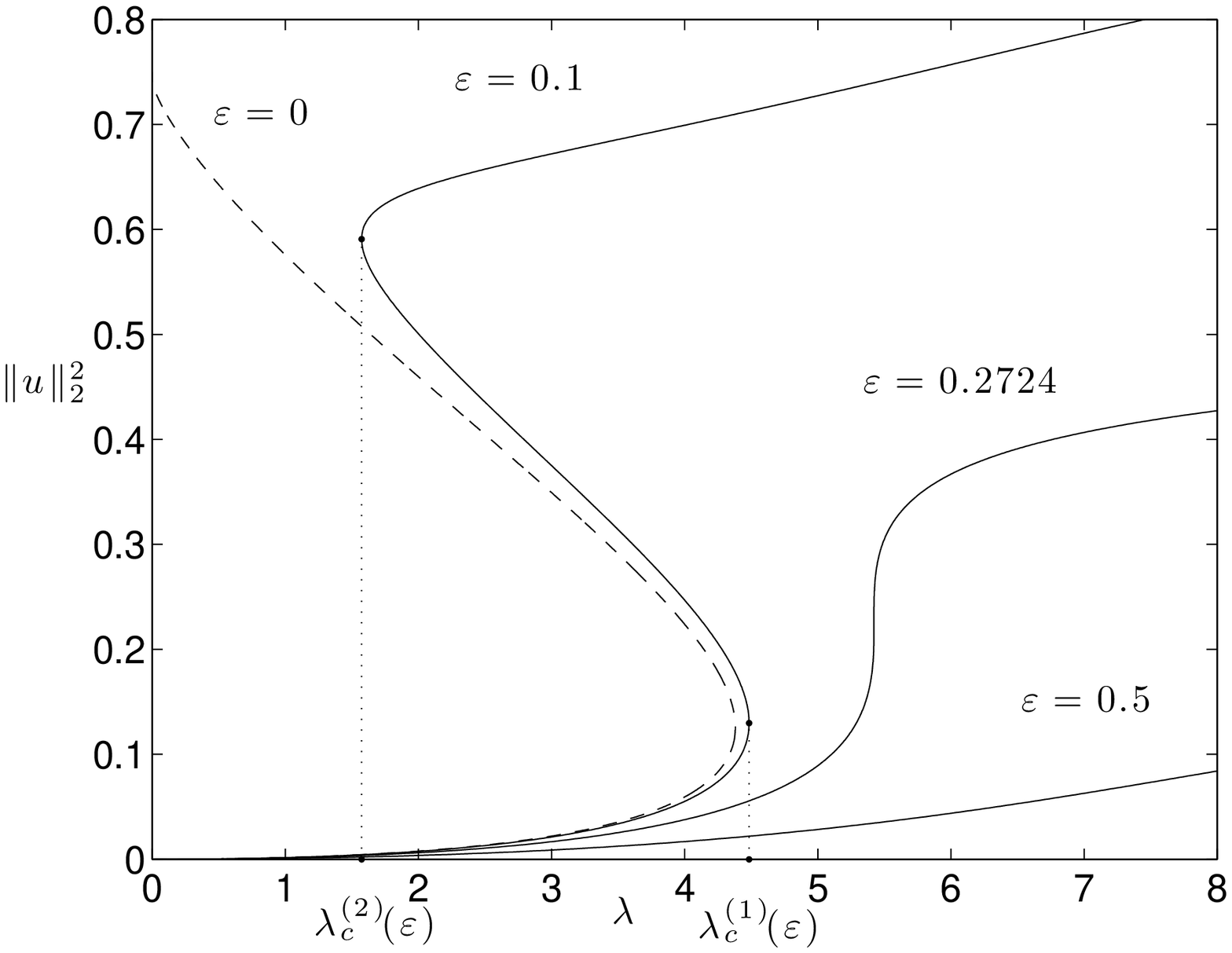}}
\parbox{5.5in}{\caption{Bifurcation \jl{diagrams showing} equilibrium solutions of \eqref{mems_intro_3} \jl{for $m=4$. Left panel: Laplacian case; right panel: b}i-Laplacian case. In each of the above, solution curves are plotted for $\eps<\eps_c$, $\eps\approx\eps_c$ and $\eps>\eps_c$ to highlight the threshold of bistability. When $\eps=0$, only two branches of solutions exist (dashed curves).  \label{fig:intro_bif_diagrams}}}
\end{figure}

%In this section, we study the equilibrium structure of the 1D regularized equations defined in %\S\ref{sec:derivation}. This leads to
\jlb{One-dimensional equilibrium solutions satisfy} the second order elliptic equation
\bsub\label{eq_1}
\begin{equation}\label{eq_1a}
 u_{xx} = \frac{\lambda}{(1+u)^2} - \frac{\lambda\eps^{m-2}}{(1+u)^m},  \quad  x\in(-1,1); \qquad u(\pm1) = 0,
\end{equation}
and its fourth order equivalent
\begin{equation}\label{eq_1b}
 - u_{xxxx} = \frac{\lambda}{(1+u)^2} - \frac{\lambda\eps^{m-2}}{(1+u)^m},  \quad  x\in(-1,1); \qquad u(\pm1) = u'(\pm1) = 0.
\end{equation}
\esub
%
%A first step is to determine the set $(\lambda,u)$ of solutions to \jl{Equations \eqref{eq_1} for} %various \jl{values of} $\eps$\jl{, and summarize the results in bifurcation diagrams}.
\jlb{Figure \ref {fig:intro_bif_diagrams} shows bifurcation diagrams obtained by numerically solving the relevant boundary value problem at fixed values of $\|u\|_2^2$. Starting from $\|u\|_2^2=0$, and $\lambda=0$, the solver identifies a value of  $\lambda$ and a solution $u(x)$ for each incremental value of the $L^2$ norm of the solution. Previously accepted solutions are used to initialize subsequent searches over a predetermined range of $\|u\|_2^2$ values.}

The bifurcation diagrams shown in Fig.~\ref{fig:intro_bif_diagrams} exhibit two remarkable deviations from the standard $\eps=0$ bifurcation diagram\jla{, displayed as a dashed curve on both panels}. The first is that for $\lambda$ arbitrarily close to $0$ \jl{and $\eps$ finite}, equations \eqref{mems_intro_3} appear to have a unique \jl{equilibrium} solution - the minimal solution branch. Secondly, \jl{there  exists a parameter range where the system exhibits bistability, and thus also possesses a stable large norm branch of equilibrium solutions. More precisely}, there \jl{is a critical value} $\eps_c$ such that for $\eps<\eps_c$, equations \eqref{mems_intro_3} are bistable over a parameter range $0 < \lambda_c^{(2)}(\eps)<\lambda<\lambda_c^{(1)}(\eps)$ while for $\eps\geq \eps_c$, a unique solution is present for each $\lambda$\jlb{, including for large values of $\lambda$}. \jl{As $\eps \rightarrow 0$, the bistable region extends towards smaller values of $\lambda$, that is $\lambda_c^{(2)} \rightarrow 0$, as is further} \ala{discussed} below and in \S\ref{sec:bistable}.

\subsection{Existence of a new branch of equilibrium solutions} \label{sec:new_sol}

\begin{figure}[htbp]
\centering
\includegraphics[height=5cm,clip]{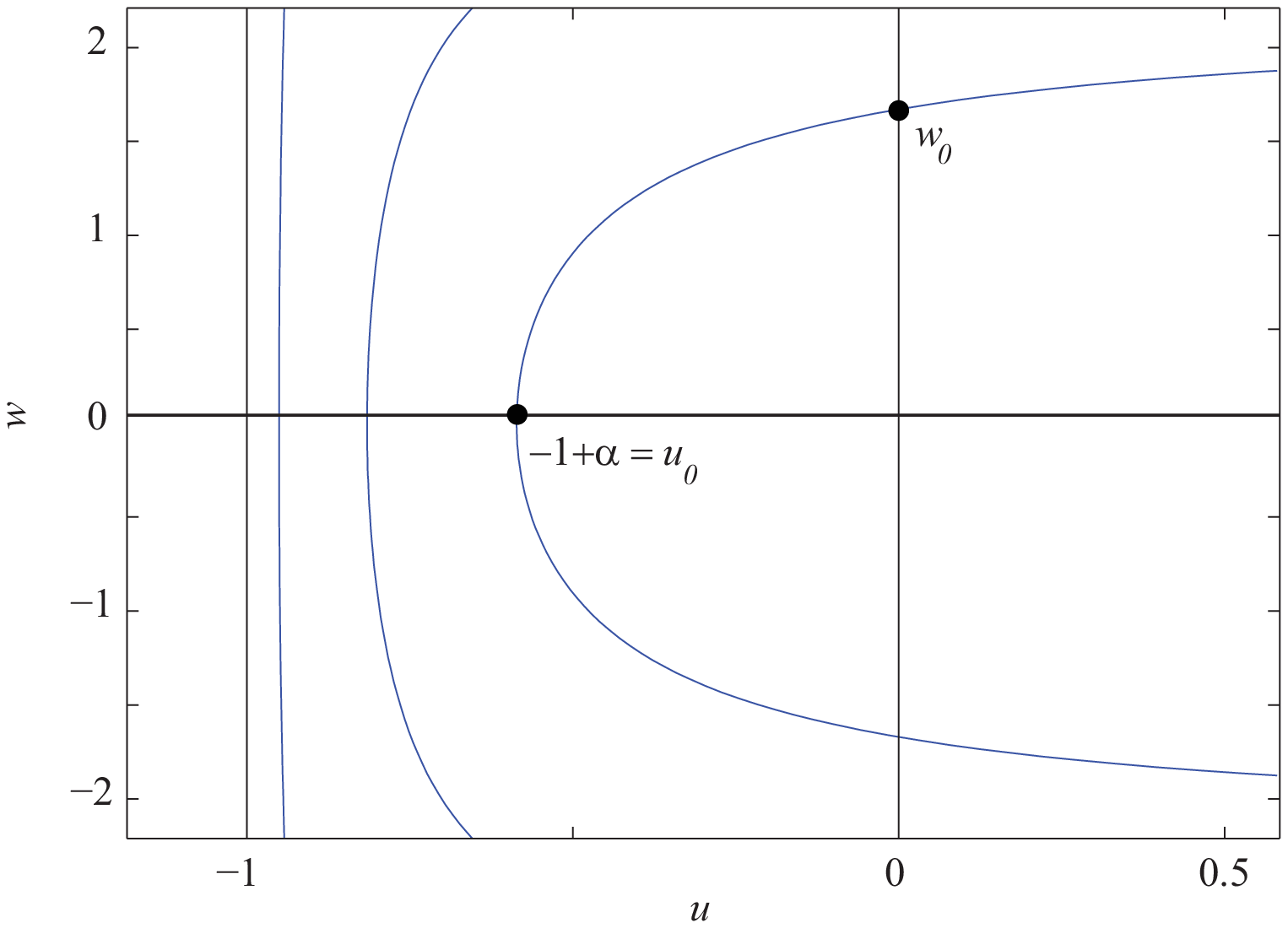} \hskip 0.8cm
\includegraphics[height=5cm,clip]{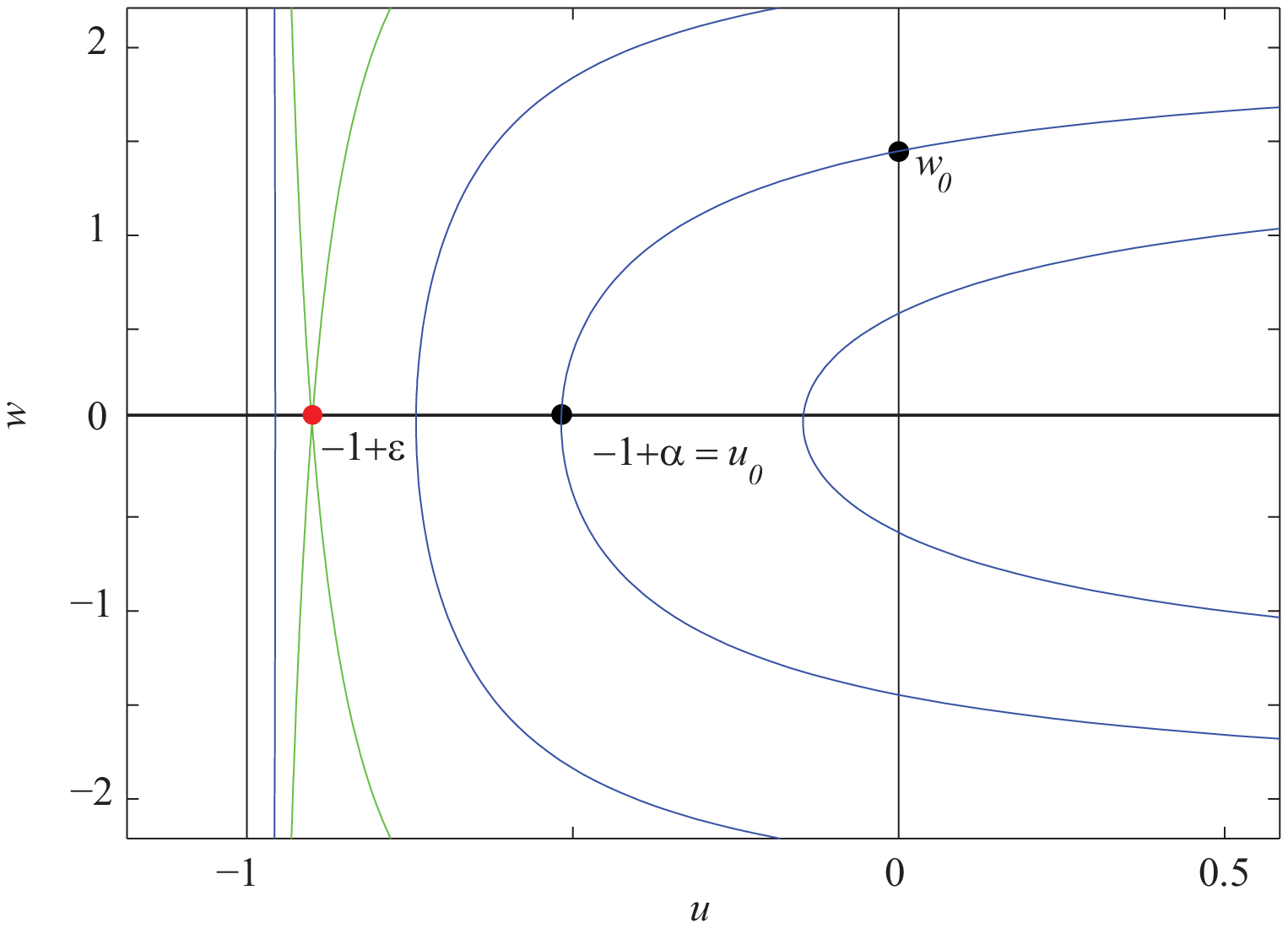}
\parbox{4.5in}{\caption{Phase portraits for the time-independent system in the Laplacian case. Left: no regularization, $\eps=0$. Right: in the presence of regularization, with $m=4$, and $\eps=0.1$. (Trajectories obtained with PPLANE)\label{fig:DS_portrait}}}
\end{figure}

To understand the existence of the saddle-node bifurcation at $\lambda= \lambda_c^{(2)}(\eps)$ when $\eps \ne 0$, we consider the dynamical system describing equilibrium solutions of Equation \eqref{mems_intro_3a}, with and without regularization. Equilibrium solutions of \eqref{mems_intro_3a} satisfy \eqref{eq_1a}, which in terms of the rescaled independent variable $y=\sqrt \lambda$ reads
\[
u_{yy} = \frac{1}{(1+u)^2}-\frac{\eps^{m-2}}{(1+u)^m}, \quad y \in [-\sqrt \lambda , \sqrt \lambda ], \quad u(\pm \sqrt \lambda) =0.
\]
The above ordinary differential equation is equivalent to the first-order system
\begin{equation}
\label{eq:1st_system}
\left \{\begin{array}{l}u_y = w \\[5pt] w_y = \ds\frac{1}{(1+u)^2}-\ds\frac{\eps^{m-2}}{(1+u)^m}\end{array}. \right.
\end{equation}
When $\eps=0$, this system has a line of singularities at $u=-1$. When $\eps \ne 0$, this line still persists, but trajectories originating near $u=0$ cannot get close to $u=-1$, due to the presence of a saddle point at $u=-1+\eps$, $w=0$ (see Figure \ref{fig:DS_portrait}). We are interested in trajectories that connect the vertical line $u=0$ to itself. Amongst these, those of half-length $\sqrt \lambda$, if any, correspond to equilibrium solutions of \eqref{eq_1a}. Note that system \eqref{eq:1st_system} is left invariant by the transformation $y \rightarrow -y$, $w \rightarrow -w$, and that the equilibrium solutions we are looking for are therefore symmetric with respect to the middle of the box. One can parameterize each trajectory that connects $u=0$ to itself by the $w$-coordinate of the point where the trajectory meets the line $u=0$ in the upper half-plane, or equivalently by the $u$-coordinate of the point where the trajectory crosses the horizontal axis. We will denote the former by $w_0$ and the latter by $u_0\equiv-1+\alpha$, with $0<\alpha \le 1$. Since distinct trajectories do not cross, $w_0$ is a decreasing function of $\alpha$ with $\alpha \in (0,1]$ for $\eps=0$ and $\alpha \in (\eps,1]$ for $\eps \ne 0$.

\begin{figure}[htbp]
\centering
\includegraphics[height=6cm,clip]{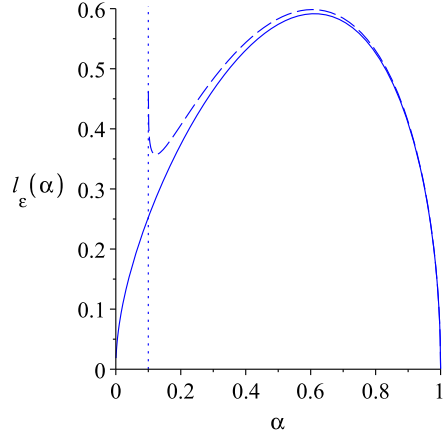}
\parbox{4.5in}{\caption{Graph of the function $l_\eps(\alpha)$ in the harmonic case in the absence of regularization ($\eps=0$, solid curve) and in the presence of regularization (for $\eps=0.1$ with $m=4$, dashed curve). The vertical line at $\alpha=\eps=0.1$ indicates where $l_\eps(\alpha)$ diverges when $\eps \ne 0$.
\label{fig:l_alpha}}}
\end{figure}

A trajectory that connects the point $(u=-1+\alpha, w=0)$ to the point $(u=0,w=w_0)$ has an equation of the form
\[
\frac{1}{2} w^2 = -\frac{1}{1+u} + \frac{\eps^{m-2}}{(m-1) (1+u)^{m-1}} + C, \qquad C = \frac{1}{\alpha} - \frac{\eps^{m-2}}{(m-1) \alpha^{m-1}},
\]
and its length $l_\eps(\alpha)$ is given by
\begin{eqnarray}
l_\eps(\alpha) &=& \int_0^{l(\alpha)} dy = \int_{-1+ \alpha }^0 \frac{du}{w} \nonumber \\
\label{eq:length} &=& \int_{-1+ \alpha }^0 \left[ \left( \frac{1}{\alpha} -\frac{1}{1+u} \right) + \frac{\eps^{m-2}}{m-1} \left( \frac{1}{(1+u)^{m-1}} - \frac{1}{\alpha^{m-1}} \right) \right]^{-1/2} du.
\end{eqnarray}
When $\eps=0$, the above integral can easily be evaluated as
\begin{eqnarray*}
l_0(\alpha) &=& \left[ \sqrt \frac{\alpha}{2} \left( \sqrt{(1+u)(1+u-\alpha)} + \alpha \ln\left(\sqrt{1+u} + \sqrt{1+u-\alpha}\right)\right)\right] _{-1+ \alpha }^0 \\
&=& \sqrt{\frac{\alpha}{2}} \left( \sqrt{1-\alpha} + \alpha \ln \left(1 + \sqrt{1-\alpha} \right) - \alpha \ln \left(\sqrt{\alpha} \right)\right).
\end{eqnarray*}
As shown in Fig.~\ref{fig:l_alpha}, for $\alpha \in (0,1]$, the graph of the above function is concave down with $ l_0(1) = 0$ and $\lim_{\alpha \rightarrow 0^+} l_0(\alpha) = 0$. It has a maximum at $\alpha_c \simeq 0.612$. As a consequence, for values of $\lambda$ such that $\sqrt{\lambda} < l_0(\alpha_c)$, there are two branches of solutions that satisfy the boundary conditions. These two branches meet when $\lambda_c = l_0(\alpha_c)^2 \simeq 0.35$. This value of $\lambda$ agrees very well with the numerically obtained value of the turning point for the bifurcation diagram of Figure \ref{fig:intro_bif_diagrams} with $\eps=0$ (dashed curve in left panel).

For $\eps \ne 0$, the change of variable $\displaystyle v= \frac{u+1-\alpha}{\alpha}$ leads to
\begin{eqnarray*}
l_\eps(\alpha) &=& \frac{\alpha^{3/2}}{\sqrt 2}\int_0^{-1+ 1/\alpha } \left[ \frac{v}{v+1} + \frac{\eps^{m-2}}{(m-1)\alpha^{m-2}} \frac{1-(1+v)^{m-1}}{(1+v)^{m-1}} \right]^{-1/2} dv \\
&=& \frac{\alpha^{3/2}}{\sqrt 2}\int_0^{-1+ 1/\alpha } \left(\frac{v}{v+1}\right)^{-1/2} \left[ 1 + \frac{\eps^{m-2}}{(m-1)\alpha^{m-2}} \frac{1-(1+v)^{m-1}}{v (1+v)^{m-2}} \right]^{-1/2} dv
\end{eqnarray*}
The above integral may be expanded in powers of $\eps$ near $\alpha =\bigoh(1)$. Since
\[
1 \le \frac{(1+v)^{m-1}-1}{v (1+v)^{m-2}} \le m-1 \qquad \hbox{for } v \ge 0,
\]
the integral appearing in the $k$-th term of the expansion is finite, and we therefore obtain a regular asymptotic expansion of $l_\eps(\alpha)$ in powers of $\eps$. For $\alpha$ near $\alpha_c$, this expansion may be used to describe how the location of the saddle node bifurcation that occurs at $\lambda = \lambda_c$ when $\eps=0$ is modified for small values of $\eps$. We indeed obtain
\begin{eqnarray*}
\lambda_c^{(1)}(\eps) &=& l_\eps(\alpha_c(\eps))^2  \\
&=& \lambda_c + \eps^{m-2}\  \frac{\alpha_c^{-m+7/2}}{m-1} \sqrt{\frac{\lambda_c}{2}} \int_0^{-1+ 1/\alpha } \left(\frac{v}{v+1}\right)^{-1/2} \frac{(1+v)^{m-1}-1}{v (1+v)^{m-2}}\, dv \\
& & + \bigoh\left(\eps^{2 (m-2)}\right),
\end{eqnarray*}
where $\alpha_c(\eps)$ is the value of $\alpha$ at which $l_\eps(\alpha)$ reaches its local maximum. For $m=4$, the above reads $\lambda_c^{(1)}(\eps) \simeq 0.350004 + 0.794451 \eps^2 + \bigoh\left(\eps^{4}\right)$, which is in agreement with the expansion of $\lambda_c^{(1)}(\eps)$ briefly mentioned in Section \ref{sec:scalings}, and derived in \cite{LLG2}.

As $\alpha \rightarrow \eps^+$, $l_\eps(\alpha)$ is expected to diverge for all values of $\eps \ne 0$, since the trajectory approaches the fixed point at $u=-1+\eps$, $w=0$. To analyze this divergence, we set $\alpha = \kappa \eps$, with $\kappa =1+\eta$ and $\eta$ small, and obtain
\[
l_\eps(\alpha) = \frac{\alpha^{3/2}}{\sqrt 2}\int_0^{-1+ 1/\alpha } \left(\frac{v}{v+1}\right)^{-1/2} \left[ g(\eta) + \frac{v\, p(v)}{(1+v)^{m-2}} \right]^{-1/2} dv,
\]
where
\begin{eqnarray*}
g(\eta) &=& \frac{1}{(1+\eta)^{m-2}} \sum_{k=1}^{m-2} {{m-2} \choose k} \eta^k\\
v\, p(v) &=& \frac{1}{(1+\eta)^{m-2}} \sum_{k=1}^{m-2} {{m-2} \choose k} \frac{k}{k+1} v^k.
\end{eqnarray*}
The function $H(v)=\displaystyle \frac{v\, p(v)}{(1+v)^{m-2}}$ is such that $H(0)=0$ and
\[
\lim_{v \rightarrow \infty} H(v)= \frac{1}{(1+\eta)^{m-2}} \frac{m-2}{m-1}.
\]
Moreover, $H$ is strictly increasing for $0 \le v \le L$, with $L=-1+1/\alpha$; a simple calculation indeed shows that its derivative is given by
\[
\frac{dH}{dv} = \frac{1}{(1+\eta)^{m-2}} \frac{1}{(1+v)^{m-1}} \left[ \frac{m-2}{2} + \sum_{k=1}^{m-3} {{m-2} \choose {k+1}} \frac{v^k}{k+2} \right] \ge \frac{m-2}{2 (1+\eta)^{m-2}}.
\]
As a consequence, on the interval $[0,L]$, $H$ is bounded above by the line tangent to its graph at the origin, and bounded below by the straight line that goes through the origin and the point of coordinates $(L,H(L))$. In other words,
\[
\frac{p(L) v}{(1+L)^{m-2} } \le H(v) \le \frac{(m-2) v}{2 (1+\eta)^{m-2}}, \qquad 0 \le v \le L.
\]
This, together with $1 \le v +1 \le L +1$ for $v \in [0,L]$, allows us to bound the term $\displaystyle \left[g(\eta) + \frac{v\, p(v)}{(1+v)^{m-2}} \right]^{-1/2}$ that appears in the expression for $l_\eps(\alpha)$, and therefore bound $l_\eps(\alpha)$. Noting that
\[
\int \frac{dv}{\sqrt{v (v + s(\eta))}} = 2 \ln\left( \sqrt v + \sqrt{v+s(\eta)}\right),
\]
we obtain $l_<(\eta) \le  l_\eps(\alpha) \le l_>(\eta)$, where $\displaystyle \eta=\frac{\alpha}{\eps}-1$ and
\begin{eqnarray*}
l_<(\eta) &=& \frac{\eps^{3/2}}{\sqrt{m-2}} \left( 1 + \frac{m+1}{2} \eta + \bigoh(\eta^2)\right) \ln\left( \frac{2(1-\eps)}{\eps \eta} + \frac{m-3}{2} \eta + \bigoh(\eta^2) \right) \\
l_>(\eta) &=& - \frac{\eps^{1/2}}{\sqrt 2} \sqrt \frac{m-1}{m-2} \ln(g(\eta)) + \bigoh\left((\eta+\eps) (\ln(\eta)+ \ln(\eps))\right).
\end{eqnarray*}
For $\eps$ fixed but small and $\eta \rightarrow 0$, we thus have
\begin{equation}
\label{eq:est_l}
- \frac{\eps^{3/2}}{\sqrt{m-2}} \ln(\eta) + \bigoh \left(\eta \ln(\eta)\right) \le l_\eps(\alpha) \le - \frac{\eps^{1/2}}{\sqrt 2} \sqrt \frac{m-1}{m-2} \ln(\eta) + \bigoh \left(\eta \ln(\eta)\right).
\end{equation}
This indicates that the graph of $l_\eps(\alpha)$ initially follows that of $l_0(\alpha)$ as $\alpha$ decreases towards $\eps$, and then diverges likes $- \ln(\eta) = -\ln(-1+\alpha/\eps)$, as shown in Fig.~\ref{fig:l_alpha}. The dashed curve is a numerical evaluation of $l_\eps(\alpha)$ for $\eps=0.1$ and $m=4$. This divergence as $\alpha \rightarrow \eps^+$ implies the existence of a third branch of solutions for $\lambda \ge \lambda^{(2)}_c(\eps)$, where $\sqrt{\lambda^{(2)}_c(\eps)}$ is the local minimum of $l_\eps(\alpha)$. The bounds in Equation \eqref{eq:est_l} show that $\sqrt{\lambda^{(2)}_c(\eps)} \rightarrow 0$ as $\eps \rightarrow 0^+$. As $\eps$ increases, the minimum of the graph of $l_\eps(\alpha)$ merges with its maximum, and only one branch of solutions exists beyond that point. This is illustrated in the numerically obtained bifurcation diagrams shown in Fig.~\ref{fig:intro_bif_diagrams} with $\eps \ne 0$ (solid curves in the left panel). The right panel of Fig.~\ref{fig:intro_bif_diagrams} shows that a similar behavior is observed in the bi-Laplacian case.

\subsection{Nature of the new branch of solutions}

\begin{figure}[H]
\centering
\subfigure[Second order, $\eps = 0.05$, $\lambda = 0.63$.]{\includegraphics[width=0.45\textwidth]{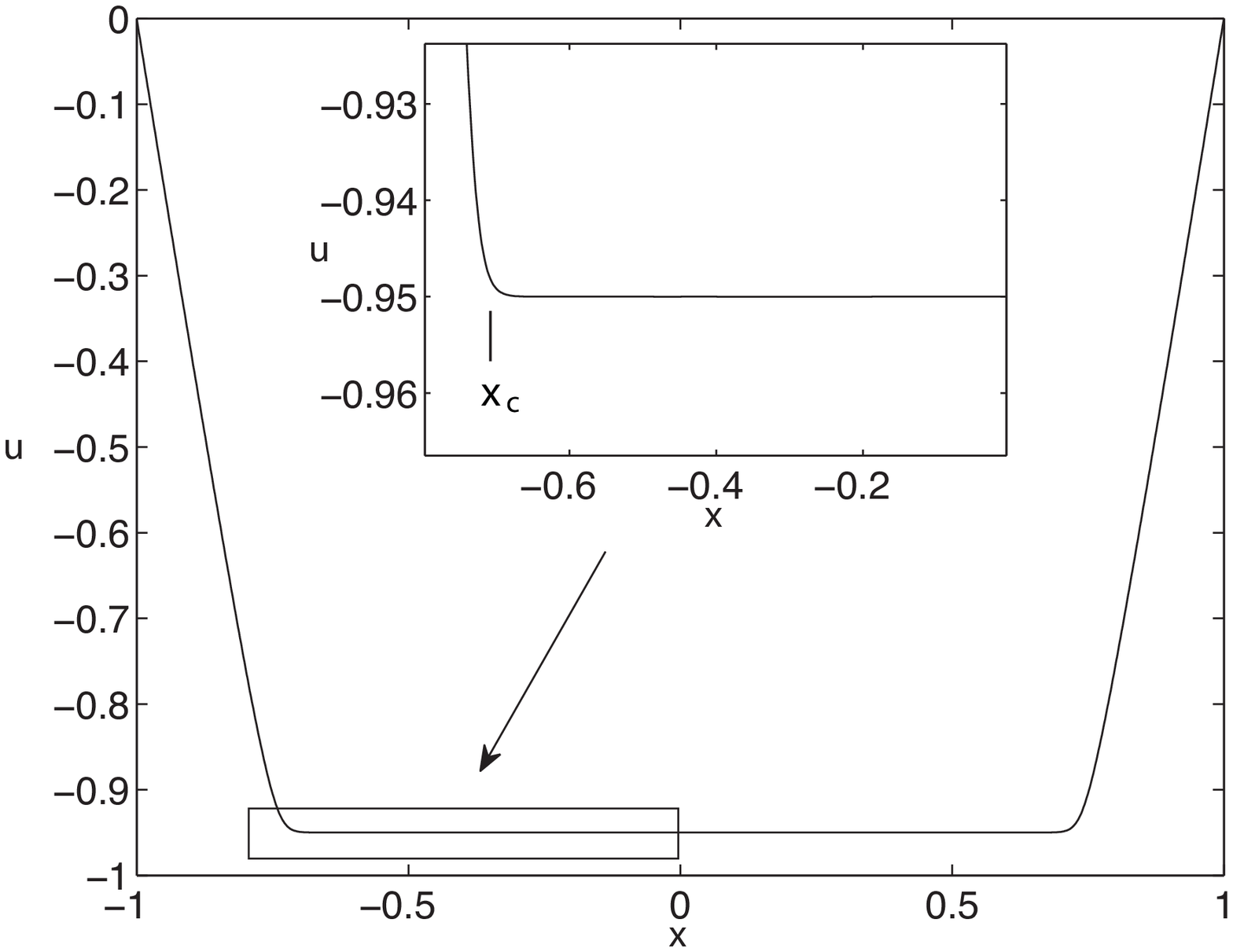}}\qquad
\subfigure[Fourth order, $\eps = 0.05$, $\lambda = 24.82$.]{\includegraphics[width=0.45\textwidth]{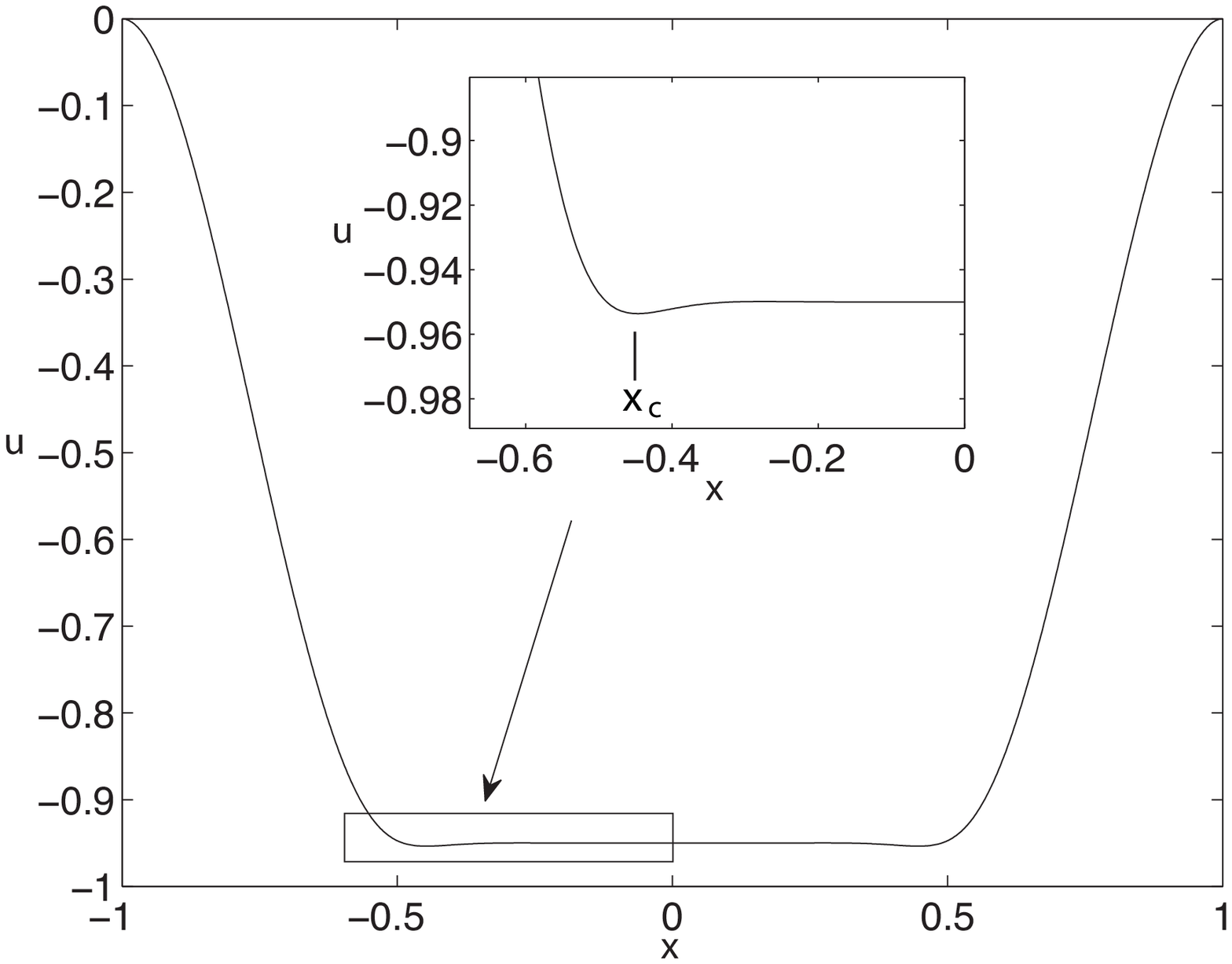}}
\parbox{5.5in}{\caption{Typical solutions of \eqref{eq_1} on the stable upper branch for $m=4$. Panels (a) and (b) represent solutions of \eqref{eq_1a} and \eqref{eq_1b} respectively. In each of the two panels, the inset panels show an enlargement of the sharp interface and touchdown region. \label{fig:soltuions_large}}}
\end{figure}

The newly present stable branch of large norm equilibria can be interpreted as a post touchdown equilibrium state. These additional solutions have three characteristic features, as illustrated in Fig.~\ref{fig:soltuions_large} \jlb{for values of $\lambda > \lambda_c$}. First, \jl{in} a large central portion \jl{of the domain}, the solution \jl{is flat and} takes on \jl{values near $-1+\eps$}. Second, a sharp transition layer \jl{links} the flat region to a profile satisfying the boundary conditions. For the Laplacian problem \eqref{eq_1a}, this sharp interface is monotone while in the bi-Laplacian case, the profile is non-monotone. Therefore, in the Laplacian case the region \jl{where $u \simeq -1 + \eps$} is spread over a finite interval while in the bi-Laplac\jl{ian} case, \jl{ $u$ attains its minimum only at} two discrete points. The third characteristic feature of \jl{this branch of equilibrium} solutions is \jl{the nature of} the profile connecting the boundary to the transition layer \jl{and in particular the size of the boundary layer.} \jlb{In what follows, we use matched asymptotic expansions to characterize these properties in the limit as $\eps \rightarrow 0$.}

\section{Scaling properties of equilibrium solutions.}\label{sec:scalings}

In this section, we construct 1D post-touchdown equilibrium \jl{configurations} of \eqref{eq_1} in the limit as $\eps\to0$. As seen in Fig.~\ref{fig:soltuions_large}, these solutions have interfaces located at $\pm x_c$, around which a narrow transition layer is centred. This transition layer separates an interior region of finite extent $(-x_c,x_c)$, from a sharp boundary profile. As explained above, the \jlb{deflection profile $u$ satisfies $u(x) = -1+ \bigoh(\eps)$ in the entire interior region, is monotonic on $[0,1]$ in the Laplacian case, and has a local minimum at the discrete points $\pm x_c$ in the bi-Laplacian case.} In both cases, it is necessary to calculate the extent of the \jl{interior region $(-x_c,x_c)$}. From numerical simulations, it \jl{appears} that $x_c$ approaches the boundary as $\eps\to0$\jl{. In the calculations below, we impose this condition, determine the scaling laws that ensue, and find the equilibrium solutions in terms of matched asymptotic expansions.}

\subsection{Laplacian Case}\label{sec:laplacian_eq}

We consider Equation \eqref{eq_1a} in the limit $\eps\to0$ and look for solutions to
\bsub\label{laplacian_eq1}
\begin{align}
\label{laplacian_eq1a}  u_{xx} = \frac{\lambda}{(1+u)^2} -  \frac{\eps^{m-2}\lambda}{(1+u)^m}, \quad x\in\jl{[-1,1]};\\[5pt]
\label{laplacian_eq1b}  u(\pm1) = 0,
\end{align}
\esub
that satisfy the following properties: (i) $u(x) = -1 + \eps + \littleoh(\eps)$ for $x \in [-x_c, x_c]$, (ii) $u(x)$ goes from its interior value of $-1 + \eps + \littleoh(\eps)$ to the value $1$ in the boundary layers $[-1,-x_c]$ and $[x_c,1]$, and (iii) there are two transition regions, centered at $\pm x_c$. From a dynamical system point of view, the particular trajectory we are interested in crosses the horizontal axis $w=u_x=0$ of the associated phase plane near but to the right of the fixed point $(-1+\eps,0)$. As the trajectory gets closer to the fixed point $(-1+\eps,0)$, the corresponding solution $u(x)$ ``spends more time'' near $u = -1 + \eps$ and therefore $x_c \rightarrow 1$. To make the scaling explicit, we write $x_c = 1- \eps^p \bar{x}_c$ for some $\bar{x}_c$ and $p$ to be determined. From symmetry considerations, since $\Omega=[-1, 1]$, we need only study the equations on the interval $[0,1]$.

\begin{figure}[H]
\centering
\includegraphics[width=0.6\textwidth]{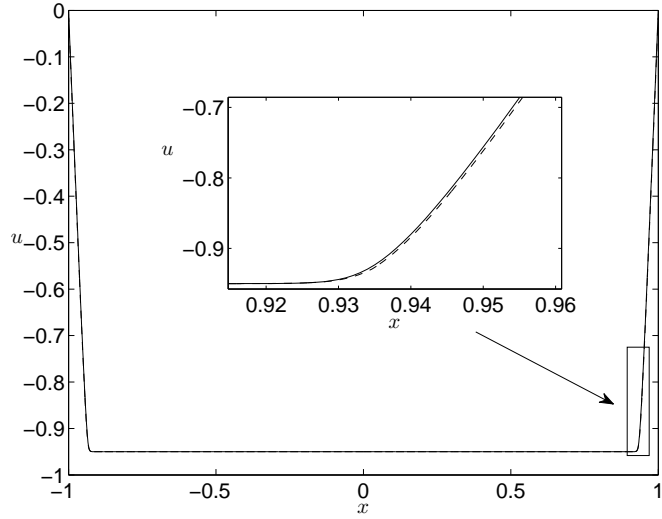}
\parbox{5.5in}{\caption{Composite asymptotic expansion of equilibrium solutions to \eqref{laplacian_eq1} for values $m=4$, $\lambda=10$, $\eps=0.05$. The solid line is the numerical solution and the dashed line is the composite asymptotic expansion. \label{Uniform_laplacian_equilibrium}}}
\end{figure}

To analyze the solution in the \jl{boundary layer} interval \jl{$[1-\eps^p \bar{x}_c,1]$}, it is convenient to use the variables
\begin{equation}\label{laplacian_eq2}
u(x) = w(\al{\eta}), \qquad \al{\eta} = \jl{\frac{x - x_c}{1 - x_c} = } \frac{x - (1-\eps^p \bar{x}_c)}{\eps^p \bar{x}_c},
\end{equation}
which transforms \eqref{laplacian_eq1} \jl{and the boundary condition $u(1-\eps^p \bar{x}_c) = -1\jlb{+ \bigoh(\eps)}$} into
\bsub
\begin{equation}\label{laplacian_eq3}
w_{\al{\eta\eta}} = \eps^{2p} \al{\lambda_c} \left[ \frac{\al{1}}{(1+w)^2} -  \frac{\eps^{m-2}}{(1+w)^m}\right], \quad \al{\eta}\in\jl{[0,1]}; \qquad w(0) = -1\jlb{+ \bigoh(\eps)}, \quad w(1) = 0,
\end{equation}
\al{where we have defined}
\begin{equation}\label{laplacian_eq3_b}
\lambda_c = \lambda \bar{x}_c^2.
\end{equation}
\esub
\jl{In light of Equation \eqref{expansion_laplacian_b} below, the introduction of $\lambda_c$ should be viewed as equivalent to expanding $x_c$ in powers of $\eps$ and $\eps \log\eps$.}
\al{We now develop the asymptotic expansion
\bsub\label{expansion_laplacian}
\begin{align}
\label{expansion_laplacian_a} w &= w_0 + \eps^{2p}\log\eps\, w_{1/2} + \eps^{2p} w_{1} + \littleoh(\eps^{2p})\\[5pt]
\label{expansion_laplacian_b} \lambda_c &= \lambda_{0c} + \eps^{2p}\log\eps\, \lambda_{1c} + \eps^{2p}\lambda_{2c} + \littleoh(\eps^{2p})
\end{align}
\esub
for solutions to \eqref{laplacian_eq3}. The $\bigoh(\eps\log\eps)$ terms are known as \emph{logarithmic switchback terms} and have previously appeared in the asymptotic construction of singular solutions to non-regularized MEMS problems \cite{LW1}. Their necessity in obtaining a consistent expansion is due to a logarithmic singularity in $w_1$ and will become apparent in the process of matching to a local solution valid in the vicinity of $\eta=0$. At leading order, the solution is given by $w_0(\al{\eta}) = -1 + \al{\eta}$ while the switchback term satisfies $w_{1/2} = a_{1/2}(\eta-1)$ where $a_{1/2}$ is a constant to be determined in the matching process. The problem for $w_1$ is
\bsub
%\label{w1_laplacian}
\begin{equation}
\label{w1_laplacian_a} w_{1\eta\eta} = \frac{\lambda_{0c}}{(1+w_0)^2}, \quad 0<\eta \leq 1; \qquad w_1(1) = 0,\\[5pt]
\end{equation}
\jlb{and its solution reads}
\begin{equation}
\label{w1_laplacian_b} w_{1} = -\lambda_{0c}\log\eta + a_1(\eta-1).
\end{equation}
\esub
}
In the transition layer near $\al{\eta}=0$, \jl{i.e. for $x \simeq x_c$}, we \jl{introduce} the local variables
\begin{equation}\label{laplacian_eq3b}
w(\al{\eta}) = -1 + \eps\jl{^\nu} v(\al{\xi}), \qquad \jlb{\xi = \frac{x-x_c}{\eps^q}},
\end{equation}
and \al{set the values $\nu=1$ \jlb{and $q=3/2$}}.
\al{This transforms equation \jla{\eqref{laplacian_eq3}} to
\begin{equation}\label{laplacian_eq4a}
v_{\al{\xi\xi}} = \eps^{2p -1} \lambda \left[ \frac{1}{v^2} - \frac{1}{v^m}\right], \quad -\infty <\al{\xi}<\infty.
\end{equation}
To balance the left and right hand sides of \eqref{laplacian_eq4a} as $\eps\to0$, the value $p=1/2$ is required.} \jlb{In order to match with the far-field solutions, we need to impose
\[
\lim_{\xi \to -\infty}v(\xi) = 1 + \littleoh(1); \qquad -1 + \eps v\left(\frac{\eta \bar x_c}{\eps}\right) \sim w(\eta) \hbox{ as } \xi = \frac{\eta \bar x_c}{\eps} \to \infty.
\]
Since the associated dynamical system has only one fixed point at ($1$,$0$) in the ($v$,$v_{\xi}$) phase plane, there is no trajectory that exactly meets these conditions. However, the unstable manifold of the above fixed point satisfies the zeroth order equation and boundary conditions. We then look for approximate solutions that solve the differential equation to a given order in $\eps$ and also have the correct behavior as $\xi \to -\infty$, to the same order in $\eps$. In particular, if the $\littleoh(\eps)$ term that appears in the boundary condition is small beyond all orders in $\eps$, we will have $v(\xi)= v_0(\xi) +  \littleoh(\eps^k)$, for all integers $k \ge 1$. The leading order problem for $v_0(\xi)$ reads
\bsub\label{laplacian_eq5}
\begin{align}
\label{laplacian_eq5a} v_{0\al{\xi\xi}} = \lambda \left[ \frac{1}{v_0^2} - \frac{1}{v_0^m}\right], \quad -\infty <\al{\xi}<\infty, \\[5pt]
\label{laplacian_eq5b} v_0(\al{\xi}) \to 1, \quad v_{0\xi}(\xi) \to 0, \quad \hbox{as } \xi \to -\infty.
\end{align}
\esub
As mentioned above, its solution corresponds to the positive branch of the unstable manifold of the fixed point ($v_0 = 1$, $v_{0\xi} = 0$) in the ($v_0$,$v_{0\xi}$) phase plane of the associated dynamical system. The above equation may be integrated once to give
\[
\frac{1}{2} v_{0\xi}^2 = \lambda \left[ - \frac{1}{v_0} + \frac{1}{(m-1)v_0^{m-1}}\right] + C_0, \qquad C_0 =  \lambda \frac{m-2}{m-1},
\]
where the value of $C_0$ was obtained from the condition as $\xi \to -\infty$. From this equation, we can infer the behavior of the unstable manifold as $\xi \to \infty$: setting $v_0(\xi) = \alpha \xi + \beta \log\xi + \bigoh(1)$ and equating the constant terms and the terms in $1/\xi$, we find
\[
v_0(\xi) = \sqrt{\frac{2 \lambda (m-2)}{m-1}} \xi - \frac{m-1}{2 (m-2)} \log\xi + \ala{\gamma + \bigoh\Big(\frac{1}{\xi}\Big)} \hbox{ as } \xi \to \infty.
\]
To match with the boundary layer expansion, we re-write $-1 + \eps v_0(\xi) + \bigoh(\eps)$ in terms of $\eta = \eps \xi / \bar x_c$ and obtain, after making use of
\[
\bar x_c = \sqrt\frac{\lambda_c}{\lambda} = \sqrt\frac{\lambda_{0c}}{\lambda} \left[1 +  \frac{\lambda_{1c}}{\lambda_{0c}} \eps \log\eps + \ala{ \frac{\lambda_{2c}}{\lambda_{0c}}\eps + \littleoh(\eps)}\right]^{1/2},
\]
the following expansion, as $\xi \to \infty$:
\begin{eqnarray*}
-1+ \eps v_0(\xi) &\simeq& -1 + \eta \sqrt{\frac{2 \lambda_{0c} (m-2)}{m-1}} + \eta \sqrt{\frac{2 \lambda_{0c} (m-2)}{m-1}} \frac{\lambda_{1c}}{2 \lambda_{0c}} \eps \log\eps + \frac{m-1}{2 (m-2)} \eps \log\eps\\
&& - \frac{m-1}{2 (m-2)} \eps \log\eta + \ala{\eps\left(\sqrt{\frac{2\lambda_{0c}(m-2)}{m-1}} \frac{\lambda_{2c}}{2\lambda_{0c}} \eta + \gamma - \frac{m-1}{4 (m-2)}\log\Big(\frac{\lambda_{0c}}{\lambda}\Big) \right)}\\
&& + \littleoh(\eps).
\end{eqnarray*}
To match with
\[
w(\eta) = -1 + \eta + \eps \log\eps\, a_{1/2} (\eta -1) - \lambda_{0c} \eps \log\eta + \eps a_1 (\eta -1) + \littleoh(\eps),
\]
we need to impose
}
\al{\begin{equation}\label{laplacian_eq7}
\lambda_{0c} = \frac{m-1}{2(m-2)}.
\end{equation}
\jlb{We then have
\[
-1+ \eps v_0(\xi) \simeq -1 + \eta + \eps \log\eps\, \Big(\eta \frac{\lambda_{1c}}{2 \lambda_{0c}} + \lambda_{0c} \Big) - \lambda_{0c} \eps \log\eta\, \ala{+ \eps\Big( \frac{\lambda_{2c}}{2\lambda_{0c}}\eta + \gamma - \frac{\lambda_{0c}}{2}\log\Big(\frac{\lambda_{0c}}{\lambda}\Big) \Big) + \littleoh(\eps)},
\]
which also requires that
\[
a_{1/2} = - \lambda_{0c}, \qquad \lambda_{1c} =  -2 \lambda_{0c}^2, \ala{\qquad \lambda_{2c} = 2a_1\lambda_{0c},\hbox{ and } a_1 = \frac{\lambda_{0c}}{2} \log\Big(\frac{\lambda_{0c}}{\lambda}\Big) - \gamma.}
\]
}
From \eqref{laplacian_eq3_b}, the two term expansion of $\bar{x}_c$ is then
\begin{equation}\label{final_contact_laplacian}
\bar{x}_c  = \left[ \frac{\lambda_{0c}}{\lambda} - 2\eps\log\eps \frac{\lambda_{0c}^2}{\lambda} + \ala{\frac{2 a_1\lambda_{0c}}{\lambda}\eps + \littleoh(\eps)} \right]^{1/2}  = \sqrt{\frac{\lambda_{0c}}{\lambda}} \Big[ 1- \lambda_{0c} \eps\log\eps + \ala{a_1\eps + \littleoh(\eps)} \Big].
\end{equation}
}
To summarize, \jlb{we expect} the equilibrium solution $u$ of \eqref{eq_1a} \jlb{to} satisfy the following properties in the limit $\eps \rightarrow 0$:
\begin{itemize}
\item $u(x) = -1 + \eps  + \littleoh(\eps^k)$, $k > 2$ in the interior region $x \in [0, x_c]$, with $\displaystyle x_c = 1 - \eps^{1/2} \bar x_c$ and $\bar x_c$ given by \eqref{final_contact_laplacian};
\item \jlb{$u(x) = -1 + \eps\, v_0(\xi) +  \littleoh(\eps)$} in the transition layer near $x_c$, with \jlb{$\displaystyle \xi = \frac{x - x_c}{\eps^{3/2}}$}.
\item \jlb{$\displaystyle u(x) = -1 + \eta - \eps \log\eps \lambda_{0c} (\eta -1) - \lambda_{0c} \eps \log\eta + \eps a_1 (\eta -1) + \littleoh(\eps)$} in the boundary layer $x \in (x_c, 1]$, \jlb{with $ \displaystyle \eta= \frac{x - x_c}{\eps^{1/2} \bar x_c}$ and $\bar x_c$ given by \eqref{final_contact_laplacian}.}
\end{itemize}
\jlb{Figure \ref{Uniform_laplacian_equilibrium} shows a comparison between the above composite asymptotic expansion and a numerical solution of the full problem, indicating very good agreement.} \jlb{In order to plot the solution obtained with matched asymptotic expansions, we have assumed that the contact point $x_c$ coincides with the maximum of the second derivative of $u(x)$, ie.,
\[
x_c = \{ x\in [0, 1] \ | \ u''(x) = \max_{y\in\Omega} u''(y) \},
\]
and calculated numerically the value of $\gamma$ in \eqref{final_contact_laplacian} accordingly.}

\begin{figure}[H]
\centering
\includegraphics[width=0.445\textwidth]{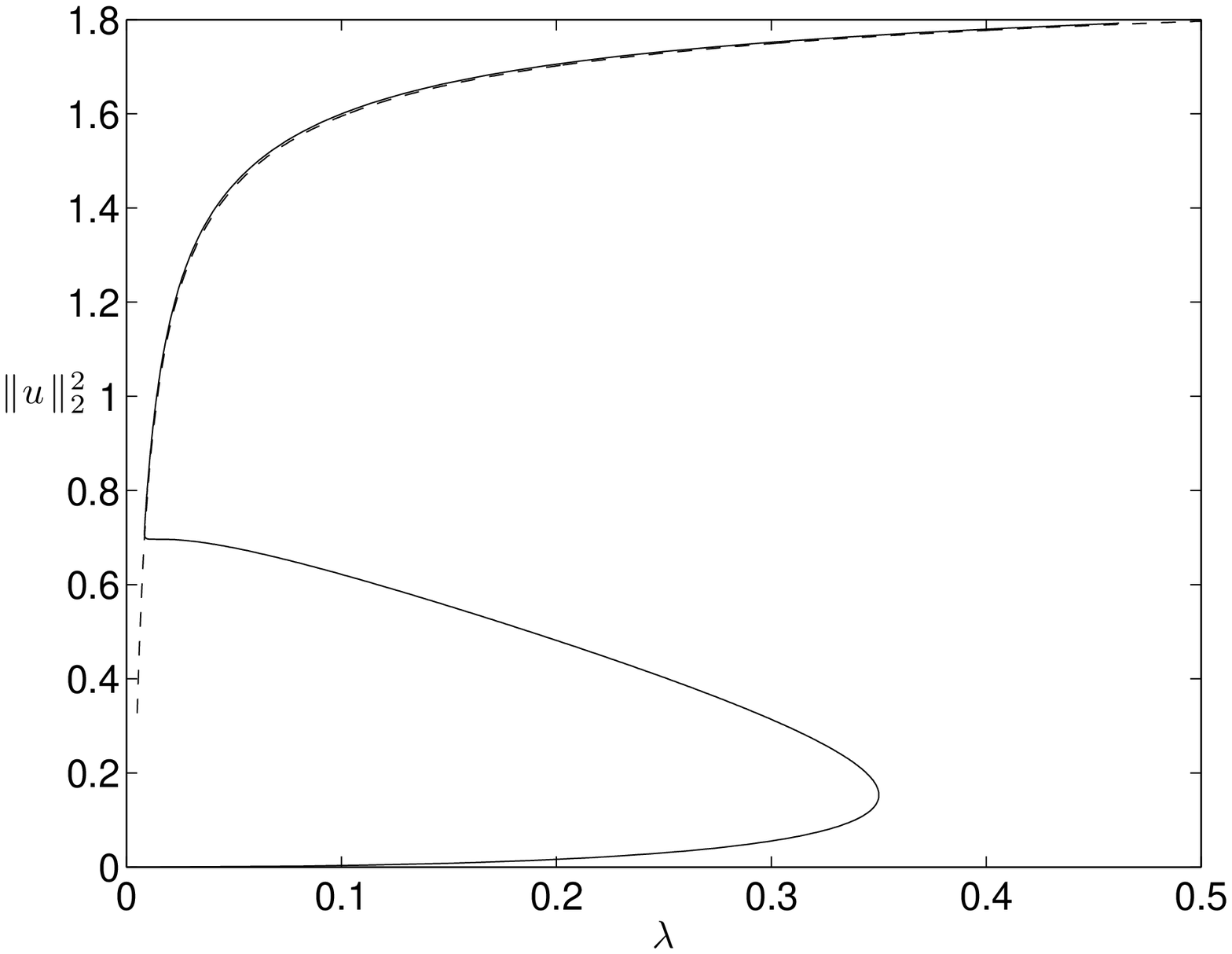} \qquad
\includegraphics[width=0.45\textwidth]{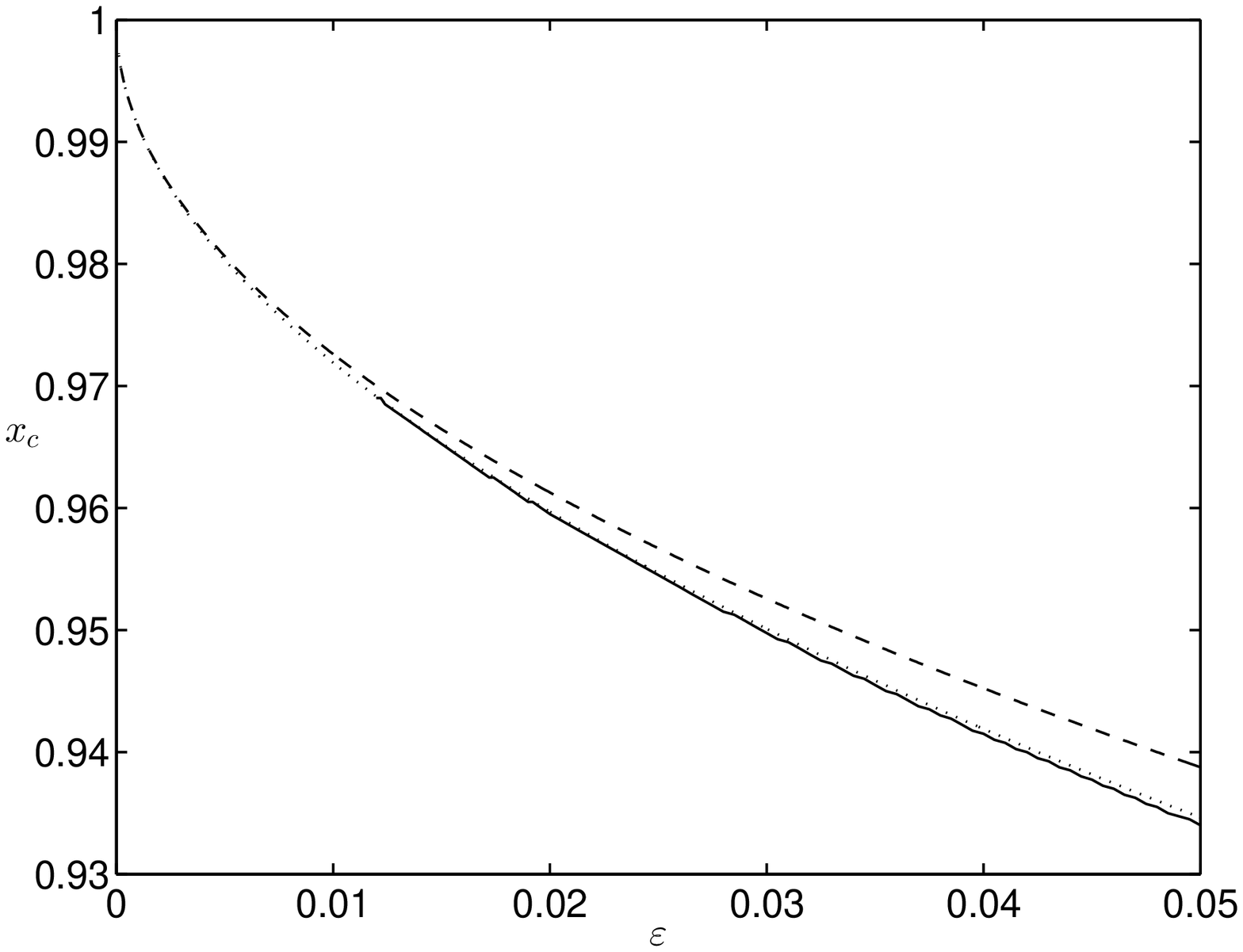}
\parbox{5.5in}{\caption{Numerical verification of \eqref{laplacian_eq8} \jl{and $p = 1/2$} for $m=4$. The left panel displays the bifurcation diagram for $\eps = 0.01$. The solid line represents the \jl{numerically obtained branches of solutions,} while the dashed line is the asymptotic formula for the large norm solution, as derived in \eqref{laplacian_eq8}. The right panel displays a comparison of the \jl{predictions for the} equilibrium contact point $x_c\jl{= 1 - \sqrt{\eps} \bar x_c}$ \jl{with $\bar x_c$ given by \eqref{final_contact_laplacian},} for fixed $\lambda = 10$ and a range of $\eps$. The dashed line is the leading order expansion while the dotted is the three term.}\label{fig:Lap_comp}}
\end{figure}

For comparison to the bifurcation diagrams, the squared $L^{\jl{2}}$ norm of \jl{the} equilibrium solution \jl{to \eqref{eq_1a} is} computed to be\jl{, in the limit $\eps \rightarrow 0$,}
\begin{align}
\nonumber \int_{-1}^1 u\jl{(x)}^2 \, dx &= 2\left[ \int_{0}^{1-\eps^{1/2}\bar{x}_c} u(x)^2 \jl{\, dx} + \int_{1-\eps^{1/2}\bar{x}_c}^1 u(x)^2 \jl{\, dx} \right] \\
\nonumber &= 2\left[\Big(-1+\eps \jl{+ \littleoh(\eps)}\Big)^2 \Big(1-\eps^{1/2} \bar{x}_c \Big) \jl{+ \bigoh(\eps^{3/2})} + \eps^{1/2} \bar{x}_c \int_{0}^{1}\left(\jlb{w}(\eta)\right)^2 \,d\al{\eta} \right]\\
\nonumber &= 2\left[ \jlb{\Big(1-2 \eps + \littleoh(\eps) \Big)} \Big(1 - \eps^{1/2}\bar{x}_c \Big) + \eps^{1/2} \bar{x}_c \jlb{\left(\frac{1}{3}-\frac{2}{3} \lambda_{0c} \eps \log\eps  + \bigoh(\eps) \right)} \jl{+ \bigoh(\eps^{3/2})} \right]\\
\nonumber&= \jlb{ 2\left[1 - 2 \eps + \littleoh(\eps) +  \eps^{1/2} \sqrt{\frac{\lambda_{0c}}{\lambda}} \Big( 1- \lambda_{0c} \eps\log\eps + \littleoh(\eps \log\eps) \Big) \left(-\frac{2}{3}-\frac{2}{3} \lambda_{0c} \eps \log\eps  + \bigoh(\eps) \right)\right]} \\
%\nonumber&= \jlb{ 2\left[1 - 2 \eps + \littleoh(\eps) +  \eps^{1/2} \sqrt{\frac{\lambda_{0c}}{\lambda}} \Big( -\frac{2}{3} + \littleoh(\eps \log(\eps)) \Big)\right]} \\
\nonumber&= \jlb{ 2\left[1 - 2 \eps + \littleoh(\eps) -\frac{2}{3} \eps^{1/2} \sqrt{\frac{\lambda_{0c}}{\lambda}} + \littleoh(\eps^{3/2} \log\eps)\right].} \\
\end{align}
If we replace \jlb{$\lambda_{0c}$} by its expression given in \jlb{\eqref{laplacian_eq7}}, the above equation reads
\begin{equation}\label{laplacian_eq8}
\|u\|_2^2 = 2 \left[1 - \frac{2\, }{3} \sqrt\frac{m-1}{2\lambda (m-2)}\, \eps^{1/2}  - \al{2\eps} + \bigoh(\eps^{3/2}\log\eps)\right].
\end{equation}
The dashed curve in the left panel of Fig.~\ref{fig:Lap_comp} shows the above quantity as a function of $\lambda$ for $m=4$ \jl{and $\eps = 0.01$}, and matches the upper branch of the bifurcation diagram very well. The right panel of Figure \ref{fig:Lap_comp} is a numerical confirmation of the $p=1/2$ scaling.

\subsection{Bi-Laplacian Case}\label{sec:bilaplacian_eq}

\jl{We now turn to} 1D equilibrium profiles of \eqref{eq_1b} in the limit $x_c\to1$ as $\eps\to0$. \jl{As in the Laplacian case, w}e write $x_c = 1 - \eps^p \bar{x}_c$ where $p$ and $\bar{x}_c$ are parameters to be determined. \al{For this particular case, a balancing argument will provide the value $p=1/4$}. \jl{We} consider the outer solution in the interval $\jlb{[}1- \eps^p \bar{x}_c,1\jlb{]}$ and employ the rescaling
\begin{equation}\label{bilaplacian_eq1}
u(x) = w(\al{\eta}), \qquad \al{\eta} = \frac{x-(1-\eps^p \bar{x}_c)}{\eps^p \bar{x}_c},
\end{equation}
which results in
\bsub\label{bilaplacian_eq2}
\begin{align}
\label{bilaplacian_eq2a} -w_{\al{\eta\eta\eta\eta}} = \eps^{4p}\al{\lambda_c} \left[ \frac{ 1}{(1+w)^2} -  \frac{\eps^{m-2}}{(1+w)^m}\right], \quad \al{\eta} \in\jlb{[0,1]}; \\[5pt]
\label{bilaplacian_eq2b} \qquad w(0) = -1, \quad w'(0) = 0, \quad w(1) = w'(1) = 0,
\end{align}
\al{where in addition, the parameter $\lambda_c$ is defined by}
\begin{equation}\label{bilaplacian_eq3c}
\lambda_{c} = \lambda \bar{x}_c^4.
\end{equation}
\esub
A logarithmic singularity also arises in the fourth order case, and as before, switchback terms are required in the expansion of \eqref{bilaplacian_eq2}. In addition there is a term at $\bigoh(\eps^{1/2})$ which arises from the translation invariance of the inner problem. In the end, the expansions
\bsub\label{bilaplacian_expansion}
\begin{eqnarray}
\label{bilaplacian_expansion_a} w &=& w_0 + \eps^{1/2} w_{1/4} + \eps^{4p}\log\eps \, w_{1/2} + \eps^{4p} \, w_1 + \littleoh(\eps^{4p}); \\
\label{bilaplacian_expansion_b} \lambda_c &=& \lambda_{0c}\, \jlb{+\, \eps^{1/2} \lambda_{1c}  + \eps^{4p} \log\eps\, \lambda_{2c} + \bigoh(\eps^{4p})} \nonumber
\end{eqnarray}
\esub
are applied to \eqref{bilaplacian_eq2}\jlb{. At leading order} $w_{0\al{\eta\eta\eta\eta}}=0$ and, with boundary conditions applied, reduces to $w_0 = -1 + 3\al{\eta}^2 - 2\al{\eta}^3$. The switchback term $w_{1/2}$ solves the problem
\bsub\label{bilap_switchback}
\begin{equation}
\label{bilap_switchback_a} w_{1/2\eta\eta\eta\eta} =0, \quad \eta \in (0,1); \qquad w_{1/2}(1)=w_{1/2\eta}(1) = 0\\[5pt]
\end{equation}
\jlb{and is given by}
\begin{equation}
\label{bilap_switchback_b} w_{1/2} (\eta) = \jlb{\alpha}_1 + \jlb{\alpha}_2 \eta - (3 \jlb{\alpha}_1+2 \jlb{\alpha}_2) \eta^2 + (2 \jlb{\alpha}_1 + \jlb{\alpha}_2)\eta^3
\end{equation}
\esub
where $ \jlb{\alpha}_1$ and $ \jlb{\alpha}_2$ are constants to be determined by matching. The term $\eps^{1/2} w_{1/4}$, not present in the Laplacian analysis of \S\ref{sec:laplacian_eq}, satisfies
\bsub\label{bilap_translation}
\begin{align}
\label{bilap_translation_a} w_{1/4\eta\eta\eta\eta} =0, \quad \eta \in (0,1); \qquad w_{1/4}(0) = w_{1/4}(1)=w_{1/4\eta}(1) = 0\\[5pt]
\label{bilap_translation_b} w_{1/4} (\eta) = \xi_0 (\eta - 2 \eta^2 + \eta^3),
\end{align}
\esub
where $\xi_0$ is a constant to be fixed in the matching procedure. The correction term at $\bigoh(\eps^{4p})$ \jlb{solves}
\begin{equation}\label{bilaplacian_correction}
-w_{1\eta\eta\eta\eta} = \frac{\lambda_{0c}}{(1+w_0)^2},\quad \eta \in (0,1); \qquad w_1(1) = w_{1\eta}(1)=0
\end{equation}
\jlb{and includes terms in $\log\eta$. The full solution is given by
\begin{eqnarray*}
w_1(\eta) &=& \Big(2 \beta_1 + \beta_2 + \frac{5}{486}\Big) \eta^3 - \Big(3 \beta_1 + 2 \beta_2 + \frac{5}{486}\Big) \eta^2 + \beta_2 \eta + \beta_1 \\
&& + \left( \frac{16}{729} \eta^3 - \frac{2}{27} \eta^2 + \frac{2}{27} \eta - \frac{1}{54} \right) \Big( \log(3-2\eta) - \log\eta\Big),
\end{eqnarray*}
where the constants $ \beta_1$ and $ \beta_2$ are arbitrary.} \jlb{In the transition layer near $x = x_c$, we define the local variables
\begin{equation}\label{bilaplacian_eq3b}
u(x) = -1 + \eps\jl{^\nu} v(\al{\xi}), \qquad \jlb{\xi = \frac{x-x_c}{\eps^q}},
\end{equation}
and set the values $\nu=1$ and $q=p+ 1/2$.
\al{This transforms equation \jla{\eqref{eq_1b}} to}
\begin{equation}
\label{bilaplacian_eq6a} -v_{\xi\xi\xi\xi} = \lambda \eps^{4p-1} \left(\frac{1}{v^2} - \frac{1}{v^m}\right), \qquad -\infty < \xi < \infty.
\end{equation}
To make this equation independent of $\eps$, we set $p=1/4$. The far-field requirements are given by
\[
\lim_{\xi \to -\infty}v(\xi) = 1 + \littleoh(1); \qquad -1 + \eps v\left(\frac{\eta \, \bar x_c}{\eps^{1/2}}\right) \sim w(\eta) \hbox{ as } \xi = \frac{\eta \, \bar x_c}{\eps^{1/2}} \to \infty.
\]
As in the Laplacian case, we will assume that the $\littleoh(\eps)$ term that appears in the far field condition as $\xi\to-\infty$ is of order $\eps^k$ with $k$ large, or that it is small beyond all orders in $\eps$, so that $v$ approximately lies on the two-dimensional unstable manifold of the fixed point $(v=1,v_{\xi}=0, v_{\xi \xi}=0, v_{\xi \xi \xi}=0)$ of the four-dimensional phase space associated to the above differential equation. We thus seek an expression for $v$ that solves  \eqref{bilaplacian_eq6a} to a given order in $\eps$ and satisfies the far-field conditions to that order as well. Equation \eqref{bilaplacian_eq6a} may be integrated once to give
\begin{equation} \label{bilaplacian_int}
- v_{\xi\xi\xi} \, v_{\xi} + \frac{1}{2} \big(v_{\xi\xi}\big)^2 +\frac{\lambda}{v}-\frac{\lambda}{(m-1) v^{m-1}}=C,
\end{equation}
where the constant of integration $\displaystyle C=\lambda \frac{m-2}{m-1}$ is determined by the value of the left-hand-side of \eqref{bilaplacian_int} at the fixed point $(v=1,v_{\xi}=0, v_{\xi \xi}=0, v_{\xi \xi \xi}=0)$. We set
\[
v(\xi)= v_0(\xi) + \eps^{1/2} v_1(\xi)+ \eps v_2(\xi)+\bigoh(\eps^{3/2}),
\]
in \eqref{bilaplacian_int} and solve the resulting equations at each order in half-integer powers of $\eps$. Since the dominant term of $w(\eta)$ as $ \eta \to 0$ is in $\eta^2$, the zeroth order solution $v_0(\xi)$ must behave like $\xi^2$ as $\xi \to \infty$. By substituting
\[
v_0 \left( \xi \right) = b_{{0}}{ \xi }^{2}+c_{{0}} \xi +d_{{0}}+\eta_{{0}}
\log\xi + \gamma_{{0}} {\frac {\log\xi }{{ \xi }
^{2}}}+ \phi_{{0}} {\frac {\log\xi }{ \xi }}+{\frac {f_{{0}}}{
\xi }}+{\frac {g_{{0}}}{{ \xi }^{2}}}+
\bigoh\Big(\frac{\log\xi} {{ \xi }^{3}}\Big)
\]
into the leading order equation and equating similar terms in $\xi$, we find
\begin{eqnarray*}
v_{{0}} \left( \xi \right) &=& b_{{0}}{ \xi }^{2}+c_{{0}} \xi +d_{{0}}+{\frac
{\lambda}{{6 \, b_{{0}}}^{2}}} \log\xi +{\frac {{\lambda}^{2} }{360\, {b_{{0}}}^{5} }}\,{
\frac {\log\xi}{{ \xi }^{2}}}+\frac{\lambda\,c_{{0}}}{12\, {b_{{0}}}^{3} } \,{\frac {1}{\xi }}\\
&&+{
\frac {\lambda\, \left( 77\,\lambda-540\,{c_{{0}}}^{2}b_{{0}}+180\,
\delta_{{3}}{b_{{0}}}^{2}+360\,{b_{{0}}}^{2}d_{{0}} \right) }{ 21600\, {b_{{0}}
}^{5}}\, {\frac {1}{ { \xi }^{2}} }} +\bigoh\Big(\frac{\log\xi} {{ \xi }^{3}}\Big),
\end{eqnarray*}
where $2 b_0^2 = C$ and $\delta_3\equiv\delta(m-3)$ is equal to $1$ if $m=3$ and to $0$ otherwise. Similar expressions for $v_1$  and $v_2$ are obtained and provided in \ref{appendix:A}. We can then evaluate $-1+\eps v(\xi)$ as $\xi \to \infty$, write the resulting expression as a function of $\eta = \sqrt \eps\, \xi / \bar x_c$, and match with the expressions for $w_i(\eta)$ found for the boundary layer expansion. Note that $x_c$, defined in \eqref{bilaplacian_eq3c}, depends on $\lambda_c$, which itself depends on $\eps$ through Equation \eqref{bilaplacian_expansion_b}. At lowest order, we obtain
\[
w_0(\eta)=-1+ b_{{0}}{\eta}^{2} \sqrt{\frac {\lambda_{{0c}}}{\lambda}}
+ a_{{1}}{\eta}^{3} \Big(\frac{\lambda_{{0c}}}{\lambda}\Big)^{3/4}
\]
which must also be equal to $-1 + 3\al{\eta}^2 - 2\al{\eta}^3$. This fixes the values of $b_0$ and $a_1$ (the coefficient of $\xi^3$ in $v_1(\xi)$) to
\[
b_{{0}}=3\, \sqrt{\frac{\lambda}{\lambda_{{0c}}}},\qquad a_{{1}}=-2\,\Big(\frac {\lambda}{\lambda_{{0c}}}\Big)^{3/4} .
\]
With $2 b_0^2 = C = \lambda (m-2)/(m-1)$, we obtain}
\begin{equation}
\label{bilaplacian_eq7c}
\lambda_{0c} = \frac{18(m-1)}{(m-2)}.
\end{equation}
\jlb{Matching the expression for $w_{1/4}$ gives
\[
a_2 = 0, \qquad c_0=\xi_0 \left(\frac{\lambda}{\lambda_{0c}}\right)^{1/4} \qquad \lambda_{1c}=-\frac{2}{3} \xi_0 \lambda_{0c},
\]
so that}
\begin{equation}\label{bilaplacian_eq13}
\al{\lambda_{1c} = -\frac{12(m-1)}{m-2} \xi_0.}
\end{equation}
\jlb{The value of $\xi_0$ will be numerically estimated to be $\xi_0 \approx -3.77$ by imposing}
\begin{equation}
\label{bilaplacian_eq6c} v_0(0) = \min_{\xi\in\mathbb{R}} v_0(\xi).
\end{equation}
\jlb{This condition removes} the translation invariance of \eqref{bilaplacian_eq6a} and \jlb{therefore uniquely specifies the contact point.
The expression for $w_{1/2}$ reads
\[
w_{{1/2}} \left( \eta \right) =\frac{3}{2} \,{\frac {\lambda_{{2c}}}
{\lambda_{{0c}}}}{ \eta}^{2} + \frac{1}{27}\,\lambda_{{0c}}\eta- \frac{1}{27}\,\lambda_{{0c}}{\eta}^{2}-{
\frac{1}{108}}\,\lambda_{{0c}}+\bigoh(\eta^3)
\]
and gives
\begin{equation}
\alpha_1 = -\frac{\lambda_{0c}}{108}, \qquad \alpha_2 = \frac{\lambda_{0c}}{27}, \qquad \lambda_{2c}=-\frac{\lambda_{0c}^2}{162}.
\end{equation}
The $w_1$ term picks up the logarithmic singularity and reads
\begin{eqnarray*}
w_{{1}} \left( \eta \right) &=& \left( \bigoh(\eta^3) + {\frac {2}{27}}\,\lambda_{{0c}}{\eta}^{2}-{\frac {2}{27}}
\,\lambda_{{0c}}\eta+{\frac {1}{54}}\,\lambda_{{0c}} \right) \log
 \left( \eta\, \sqrt [4]{\frac{\lambda_{{0c}}}{\lambda}}
 \right) \\
&&+\left(-\frac{3}{2} \frac{\lambda_{3c}}{\lambda_{0c}}+ \frac{1}{12}\xi_0^2\right) \eta^3 + \left( -{\frac {7}{81}}\,
\lambda_{{0c}}- c_{{1}} \sqrt [4]{\frac{\lambda_{{0c}}}{\lambda}}+\frac{3}{2}\,{\frac {\lambda_{{3c}}}{\lambda_{{0c}}}} \right) {
\eta}^{2}\\
&&+ \left( -\frac{1}{6}\,{\xi_{{0}}}^{2}+ c_{{1}} \sqrt [4]{\frac{\lambda_{{0c}}}{\lambda}} \right) \eta+d_{{0}},
\end{eqnarray*}
leading to
\[
d_0= \beta_{{1}}={\frac {7}{81}}\,\lambda_{{0c}}+\frac{1}{12}\,{\xi_{{0}}
}^{2},\qquad \beta_{{2}}=-\frac{1}{6}\,{\xi_{{0}}}^{2}+\sqrt [4]{\frac{\lambda_{{0c}}}{\lambda}}\,c_{{1}},
\]
and
\[
 \lambda_{{3c}}=-{\frac {28}{243}}\,
{\lambda_{{0c}}}^{2}+\frac{1}{18}\,\lambda_{{0c}}\,{\xi_{{0}}}^{2}-{\frac {5}{
729}}\,\lambda_{{0c}}-2/3\,\sqrt [4]{\frac{\lambda_{{0c}}}{\lambda}}\,\lambda_{{0c}}\, c_{{1}}.
\]
}
Combining \eqref{bilaplacian_eq13}, \eqref{bilaplacian_eq7c}, and \eqref{bilaplacian_eq3c}, the two term expansions for the contact points are
\[
x_c = \pm \left[1 - \left[ \frac{18(m-1)}{\lambda(m-2)} \right]^{1/4}  \left( \eps^{1/4} - \frac{\xi_0}{6} \eps^{3/4} - \frac{\lambda_{0c}}{648}\eps^{5/4}\log\eps  + \bigoh(\eps^{5/4}) \right) \right].
\]
where $\xi_0 \approx -3.77$.

To summarize, \jlb{we expect} the equilibrium solution $u$ of \eqref{eq_1b} \jlb{to} satisfy the following properties in the limit $\eps \rightarrow 0$:
\begin{itemize}
\item $u(x) = -1 + \eps \jlb{ + \littleoh(\eps^k)}$, $k > 2$ in the interior region $x \in [0, x_c]$, with $\displaystyle x_c = 1 - \eps^{1/4} \bar x_c$;
\item \jlb{$u(x) = -1 + \eps\, v_0(\xi) + \eps^{3/2} v_1(\xi) +  \eps^{2} v_2(\xi) + \bigoh(\eps^{5/2})$} in the transition layer near $x_c$, with \jlb{$\displaystyle \xi = \frac{x - x_c}{\eps^{3/4}}$}.
\item \jlb{$\displaystyle u(x) = -1 + 3 \eta^2 - 2 \eta^3 + \xi_0\, \eta\, (\eta -1)^2 \eps^{1/2}+\bigoh(\eps \log\eps)$} in the boundary layer $x \in (x_c, 1]$, \jlb{with $ \displaystyle \eta= \frac{x - x_c}{\eps^{1/4} \bar x_c}$.}
\end{itemize}
For comparison \jl{with} the numerical bifurcation diagram, the squared $L^2$ norm of the \jlb{composite asymptotic expansion} is calculated to be
\begin{align}
\nonumber \|u\|^2_2 &= 2\left[ \int_{0}^{1-\eps^{1/4}\bar{x}_c} u(x)^2 \jl{\, dx} + \int_{1-\eps^{1/4}\bar{x}_c}^1 u(x)^2 \jl{\, dx} \right] \\
\nonumber &= 2\left[(-1+\eps \jl{+ \littleoh(\eps)})^2 (1-\eps^{1/4} \bar{x}_c) \jl{+ \bigoh(\eps^{\jlb{3}/4})} + \eps^{1/4} \bar{x}_c \int_{0}^{1}\left(w_0\jl{+ \eps^{1/2} w_{1/4} + \bigoh(\eps^{5/4}\log\eps)}\right)^2 \,d\al{\eta} \right]\\
\nonumber &= 2\left[ (-1+\eps\jl{+ \littleoh(\eps)})^2\Big(1 -\eps^{1/4} \bar{x}_c \Big) + \eps^{1/4} \bar{x}_c \int_{0}^{1}w_0^2\jl{+ 2\eps^{1/2} w_0w_{1/4} \,d\al{\eta} } \jl{+ \bigoh(\eps^{3/4})} \right]
\end{align}
To simplify this expression, we calculate that
\[
\int_0^1 w_0^2\,d\eta = \frac{13}{35},\qquad \int_0^1 w_0 w_{1/4}\,d\eta = -\frac{11\xi_0}{210},
\]
and apply the expansion
\begin{equation}\label{biharm_contact_two_term}
\bar{x}_c = \left[\frac{18(m-1)}{\lambda(m-2)} \right]^{1/4}  \left( 1 - \frac{\xi_0}{6} \eps^{1/2} + \bigoh(\eps\log\eps) \right),
\end{equation}
which finally results in the value
\begin{equation}\label{bilaplacian_eq4}
\|u(x;\eps)\|^2_2 = 2\left[1 - \, \frac{22}{35} \left(\frac{18(m-1)}{\lambda(m-2)}\right)^{1/4} \eps^{1/4} \jl{+ \bigoh(\eps^{3/4})}\right].
\end{equation}

\begin{figure}[H]
\centering
\includegraphics[width=0.43\textwidth]{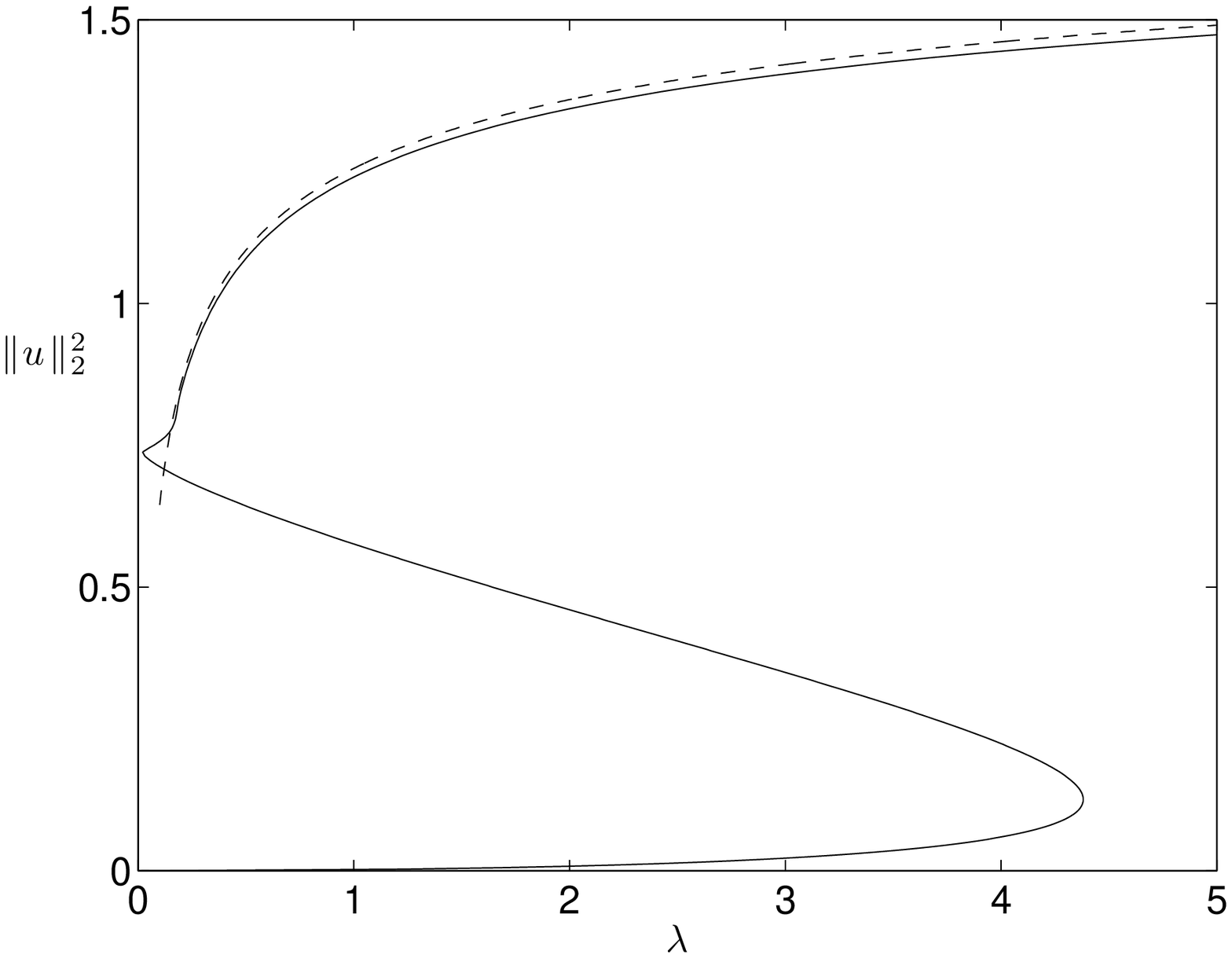} \qquad
\includegraphics[width=0.45\textwidth]{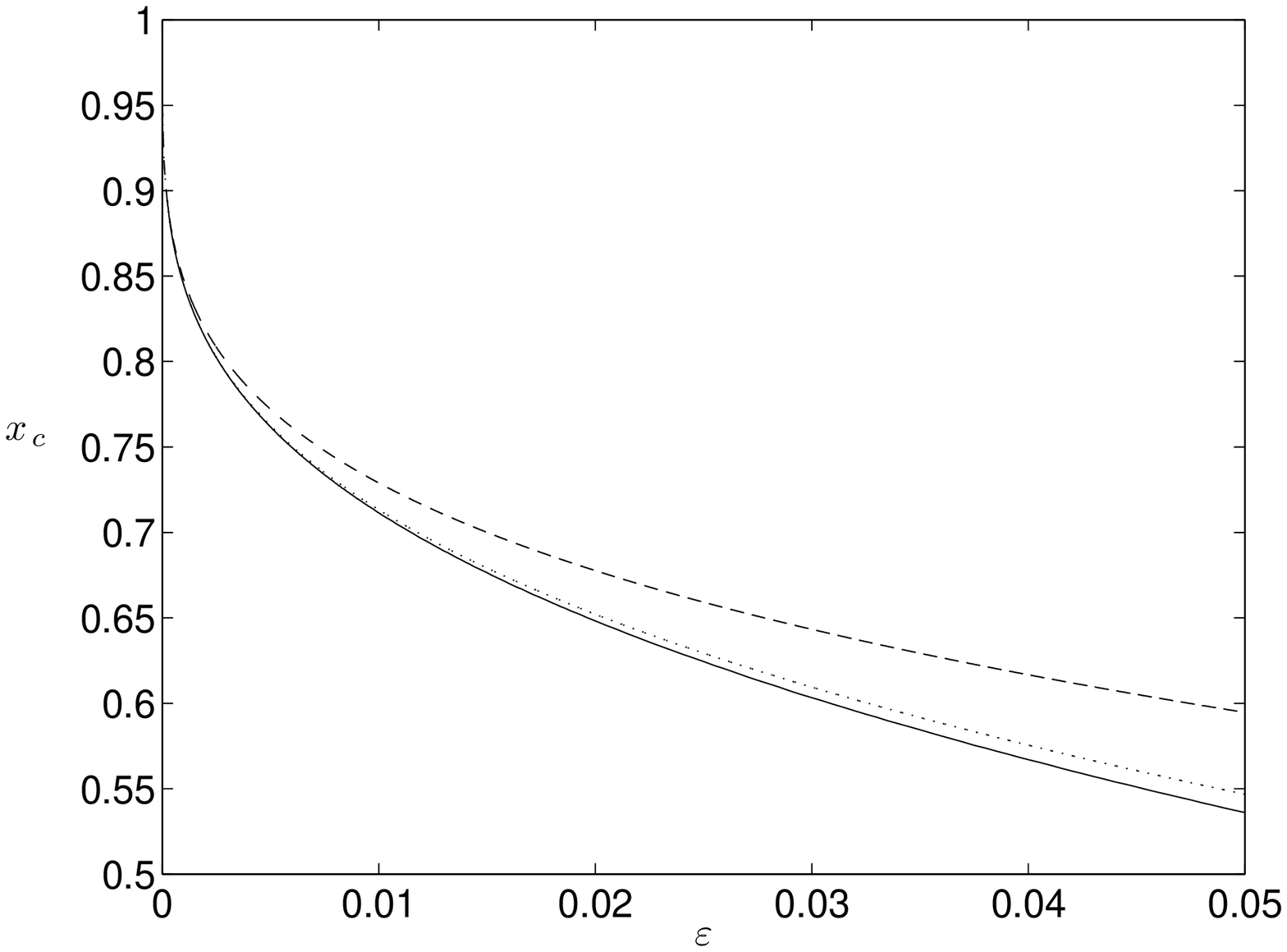}
\parbox{5.7in}{\caption{Numerical verification of asymptotic calculations for the bi-Laplacian case and $m=4$. The left panel displays the bifurcation diagram for $\eps = 0.005$. The solid line represents the \jl{numerically obtained branches of solutions,} while the dashed line is the asymptotic formula for the large norm solution, as derived in \eqref{bilaplacian_eq4}. The right panel displays a comparison of the \al{one term (dashed line) and two term (dotted line)} \jl{predictions for the} equilibrium contact point $x_c= 1 - \eps^{1/4}\bar{x}_c$ \jl{with $\bar x_c$ given by \eqref{biharm_contact_two_term},} for fixed $\lambda = 50$ and a range of $\eps$.}\label{fig:bilap_comp}}
\end{figure}

\jl{This quantity is plotted (dashed curve) in the left panel of Fig.~\ref{fig:bilap_comp} as a function of $\lambda$ for $m=4$, and is in good agreement with the numerically computed bifurcation diagram of Fig.~\ref{fig:intro_bif_diagrams}. As before, the right panel of Fig.~\ref{fig:bilap_comp} is a numerical confirmation of the $\eps$-scaling (with an exponent $p=1/4$ in this case) of the width of the boundary layer.}

\begin{figure}
\centering
\includegraphics[width=0.4\textwidth]{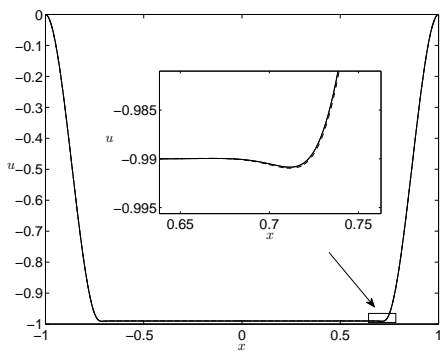}
\parbox{5.5in}{\caption{\jlb{Composite} asymptotic expansion of equilibrium solutions to \eqref{eq_1} for values $m=4$, $\lambda=50$, $\eps=0.01$. The solid line is the numerical solution and the dashed line is the composite asymptotic expansion. \label{Uniform_bilaplacian_equilibrium}}}
\end{figure}

\subsection{Singular asymptotics and bistability}\label{sec:bistable}

In this section, we briefly focus on another of the remarkable departures \jl{from} the standard $\eps=0$ bifurcation diagram \jl{displayed by} the regularized equations \eqref{eq_1}\jl{, namely} the presence of bi-stability for a certain range of $\eps$. \jl{Recall that t}he three characteristic bifurcation diagrams shown in Fig.~\ref{fig:intro_bif_diagrams} \jl{have} the following \jl{features.} For $\eps\in(0,\eps_c)$, the bifurcation diagrams of \eqref{eq_1} have two fold points $\lambda_c^{(1)}$ and $\lambda_c^{(2)}$, which results in bistable behaviour for $\lambda_c^{(2)}<\lambda<\lambda_c^{(1)}$. At the critical value $\eps=\eps_c$, there is a single cubic fold point, while for $\eps_c<\eps$, there are no fold points and \eqref{eq_1} has a unique solution for each $\lambda$.

\begin{figure}[H]
\centering
\includegraphics[width=0.45\textwidth]{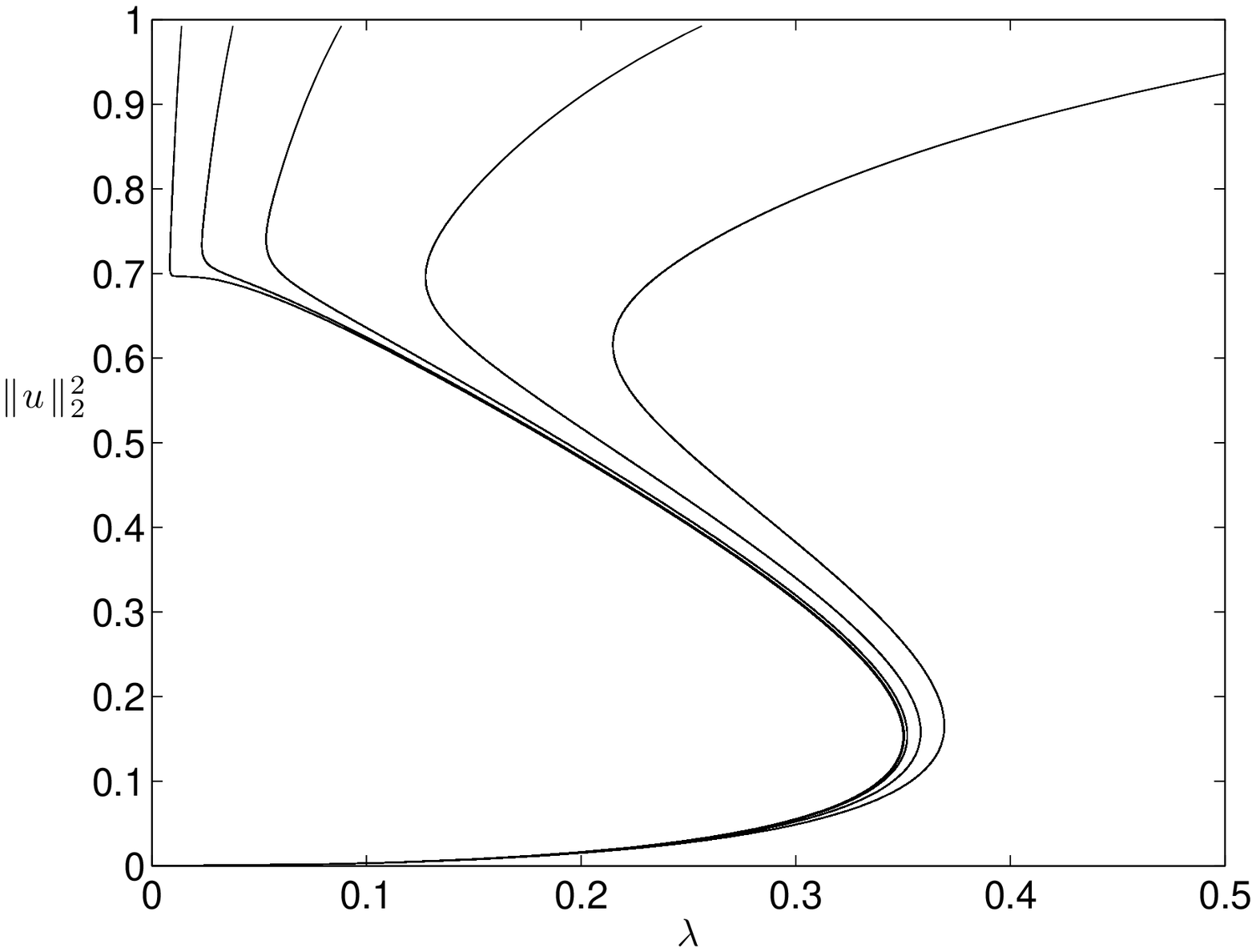} \qquad
\includegraphics[width=0.45\textwidth]{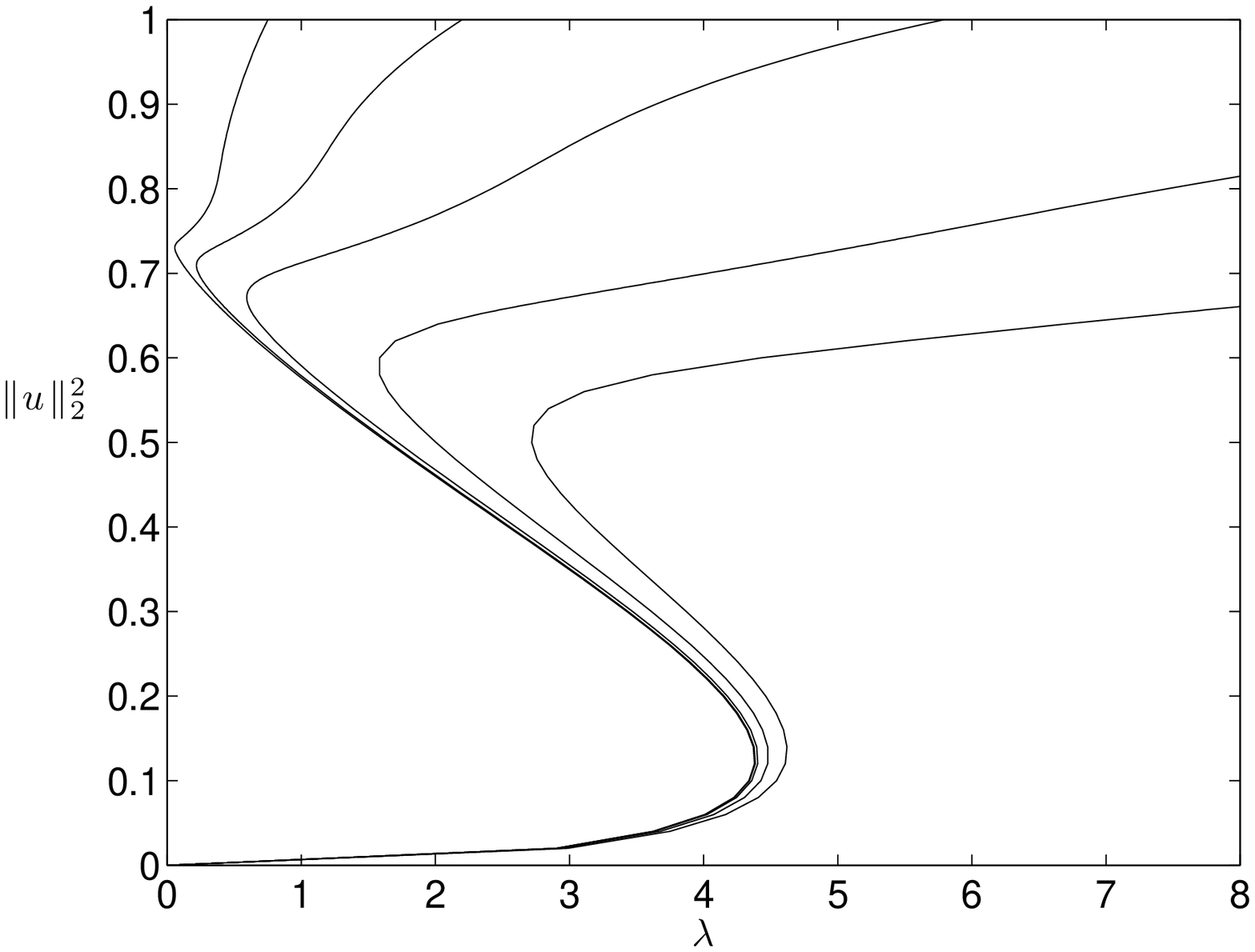}
\parbox{5.5in}{\caption{Bifurcation diagrams of \eqref{eq_1} for $\eps = 0.01,0.025, 0.05, 0.1, 0.15$ \jl{(from left to right)} and $m=4$. Left panel: Laplacian case \eqref{eq_1a}; right panel: bi-Laplacian case \eqref{eq_1b}. \label{fig:big_rang_of_eps} }}
\end{figure}

In Fig.~\ref{fig:big_rang_of_eps}, the bifurcation diagrams of \eqref{eq_1} are displayed for a range of $\eps\in(0,\eps_c)$ and $m=4$. In each case, the fold point $\lambda_c^{(2)}(\eps)$ is observed to depend quite sensitively on the parameter $\eps$, while the principal fold point $\lambda_c^{(1)}(\eps)$ exhibits smaller variations as $\eps$ increases. \ala{In essence, the regularizing term of the governing equations generates a regular perturbation to solutions of the $\eps=0$ problem whenever $1+u = \bigoh(1)$, and a singular perturbation to solutions of the $\eps=0$ problems whenever $u+1\simeq \eps$.} In each of the cases \jl{represented} in Fig.~\ref{fig:big_rang_of_eps}, the two fold points are empirically seen to be increasing functions of $\eps$\jl{, with $\lambda_c^{(2)}(\eps)$ increasing faster than $\lambda_c^{(1)}(\eps)$}. We therefore \jl{expect} the two fold points \jl{to} eventually merge at some critical $\eps_c$, where the condition
\begin{equation}\label{bistable:condition}
\lambda_c^{(1)}(\eps_c) =\lambda_c^{(2)}(\eps_c)
\end{equation}
is satisfied. \ala{The bistable features of the regularized system are very interesting as they give the device the capacity to switch robustly between a large and a small norm equilibrium state. The relative magnitude of the switching voltage required to transition the device between these two states is given, for $\eps<\eps_c$, by the quantity $\lambda_c^{(1)}(\eps)-\lambda_c^{(2)}(\eps)$.}

It is therefore desirable to obtain explicit formulae for $\lambda_c^{(1)}(\eps)$ and $\lambda_c^{(2)}(\eps)$ so that the critical parameter $\eps_c$ \jl{may} be estimated from the condition \eqref{bistable:condition} and the bistable nature of the regularized system understood. In a forthcoming paper \cite{LLG2}, a detailed singular perturbation analysis is employed to accurately locate these fold points. The main results are explicit expansions of form
\bsub\label{sing_results}
\begin{equation}\label{sing_results_a}
\lambda_c^{(1)}(\eps) \sim \lambda_0\jl{^{(1)}} + \eps^{m-2} \lambda_1\jl{^{(1)}} + \cdots, \\[5pt]
\end{equation}
for the principal fold point in the Laplacian or bi-Laplacian case. The scaling of the second fold point is quite different for the second and fourth order problems, \jl{namely}
\begin{equation}\label{sing_results_b}
\begin{array}{rlc}
\lambda_c^{(2)}(\eps) &\sim\ \lambda_0\jl{^{(2)}} \eps + \lambda_1\jl{^{(2)}} \eps^2\log\eps + \lambda_2\jl{^{(2)}} \eps^2 + \cdots &(\mbox{Laplacian})  \\[5pt]
\lambda_c^{(2)}(\eps) &\sim\ \lambda_0\jl{^{(2)}} \eps^{3/2} + \lambda_1\jl{^{(2)}} \eps^2 + \cdots & (\mbox{bi-Laplacian})
\end{array}
\end{equation}
\esub
In the above formulations, closed form expressions for the coefficients $\lambda_0\jl{^{(i)}}$, $\lambda_1\jl{^{(i)}}$ and $\lambda_2\jl{^{(i)}}$ \jl{are} established \cite{LLG2}.

\section{Discussion}\label{sec:discussion}

In this work we have proposed and analyzed a formulation for regularization of touchdown in MEMS capacitors. These considerations have resulted in a new family of models whose solutions remain globally bounded in time for all parameter regimes, followed by equilibration to new steady states. Interestingly, the presence of these new stable equilibri\jl{a} results in bistable behaviour for a range of parameter values. This may be useful in practical applications since bistable systems can be used to create robust switches. We have described how equilibrium solutions depend on the parameters $\lambda$ and $\eps$ in terms of bifurcation diagrams, for both the Laplacian and the bi-Laplacian cases. Using asymptotic analysis, we have also given a complete characterization of the scaling properties of the upper branch of equilibrium solutions, which correspond to attracting post-touchdown configurations of the regularized equations.

There are several avenues of future exploration emanating from this study. The method of regularization used in the present work is a first attempt at understanding behavior of MEMS after touchdown. It is natural to ask whether this bistability feature is generic to a larger family of regularized models.

An interesting problem is the characterization of the intermediate dynamics between the initial regularized touchdown event and the equilibration to the post touchdown states. As is typical with such obstacle type regularizations, the equations \eqref{mems_intro_3} give rise to a free boundary problem for the extent of the touchdown region, which is amenable to analysis (cf. \ref{mems_intro_3}). \jl{In a forthcoming paper \cite{LLG}, we describe the dynamic evolution of the periphery of the growing post-touchdown region, in both one and two spatial dimensions.}

\section*{Acknowledgments}
K.G. acknowledges support from NSF award DMS-0405596. A.E.L acknowledges support from the Carnegie Trust for the Universities of Scotland.

\begin{appendix}
\section{Expressions for $v_1$ and $v_2$}\label{appendix:A}

We give below the expressions for $v_1(\xi)$ and $v_2(\xi)$ such that $v = v_0 + \eps^{1/2} v_1 + \eps v_2 +\bigoh(\eps^{3/2})$ solves \eqref{bilaplacian_int} to order $\eps^{1/2}$ and $\eps$ respectively.

\begin{eqnarray*}
v_{{1}} \left( x \right) &=&a_{{1}}{\xi}^{3}+{\frac {3 a_{{1}}c_{{0}}}{2 b_{{0}}}} { \xi}^{2} +c_{{1}} \xi +d_{{1}}+{\frac {\lambda\,c_{{0}}a_{{1}}
 }{2\, {b_{{0}}}^{4}}} \log  \left( \xi \right) +{\frac {\lambda\,a_{{1}} }{{b_{{0}}}^{3}}} \xi \log \left( \xi \right) +{\frac {{\lambda}^{2}a_{{1}}}{24\, {b_{{0}}}^{6}}}
\frac{ \log  \left( \xi \right) }{ \xi }\\
&& + \gamma_{{1}} {\frac {\log \left( \xi \right) }{{ \xi }^{2}}}-{\frac {\lambda\, \left( -36\,{c_{{0}}}^{2}a_{{1}}b_{{0}}+72\,d_{{0}}a_{{1}}{b_{{0}}}^{
2}-24\,c_{{1}}{b_{{0}}}^{3}-25\,\lambda\,a_{{1}}+36\,\delta_{{3}}a_{{1}}{b_{{0}}}^{2} \right) }{ 288\, {b_{{0}}}^{6} \xi }}\\
&& +{\frac {g_{{1}}}{{ \xi }^{2}}} +\bigoh\Big({\frac {\log \left( \xi \right) }{{ \xi }^{3}}}\Big),
\end{eqnarray*}
and
\begin{eqnarray*}
 v_{{2}} \left( \xi \right) &=& a_{{2}}{ \xi }^{3}+b_{{2}}{ \xi }^{2}+c_{{2}} \xi +d_{{2}}+\kappa_{{2}} \left( \log  \left( \xi \right)  \right) ^{2}+\eta_{{3}}\,
\log  \left( \xi \right) +\eta_{{4}} \, \xi \log  \left( \xi \right) +\eta_{{5}}\, { \xi }^{2}\log  \left( \xi \right)\\
&& + \phi_{{2}} {\frac {\log  \left( \xi \right) }{ \xi}}+ \gamma_{{2}} {\frac {\log  \left( \xi \right) }{{ \xi }^{2}}}+{\frac {f_{{2}}}{ \xi }}+{\frac {g_{{2}}}{{ \xi }^{2}}}+ \bigoh\Big({\frac {\log \left( \xi \right) }{{ \xi }^{3}}}\Big),
\end{eqnarray*}
where
\begin{eqnarray*}
\eta_{{3}}&=& \lambda\,{\frac {\left( -18\,\delta_{{3}}{a_{{1}}}^{2}{b_{{0}}}^{2}+16\,\lambda\,{a_{{1}}}^{2}+9\,a_{{2}}c_{{0}}{b_{{0}}}^{3}+9\,{a_{{1}}}^{2}{c_{{0}}}^{2}b_{{0}}-36\,{a_{{1}}}^{2}{b_{{0}}}^{2}d_{{0}}+9\,a_{{1}}c_{{1}}{b_{{0}}}^{3} \right) }{18\, {b_{{0}}}^{7}}},\\
\eta_{{4}}&=&{\frac {6\,\lambda\,c_{{0}}{a_{{1}}}^{2}+4\,\lambda\,a_{{2}}{b_{{0}}}^{2}}{4 \, {b_{{0}}}^{5}}},\qquad
\eta_{{5}}={\frac {3 \lambda\,{a_{{1}}}^{2}}{2 \, {b_{{0}}}^{4}}},\qquad
\kappa_{{2}}={\frac {{\lambda}^{2}{a_{{1}}}^{2}}{12 \, {b_{{0}}}^{7}}},\\
b_{{2}}&=&-{\frac {14\,\lambda\,{a_{{1}}}^{2}-12\,a_{{2}}c_{{0}}{b_{{0}}}^{3}+9\,{a_{{1}}}^{2}{c_{{0}}}^{2}b_{{0}}-12\,a_{{1}}c_{{1}}{b_{{0}}}^{3}}{8 \, {b_{{0}}}^{4}}},\\
\phi_{{2}}&=&-{\frac {-4\,{\lambda}^{2}a_{{2}}{b_{{0}}}^{2}+7\,{\lambda}^{2}c_{{0}}{a_{{1}}}^{2}+720\ a_{{1}}\gamma_{{1}}{b_{{0}}}^{7}} {96 \, {b_{{0}}}^{8}}},\\
f_{{2}}&=&{\frac {\lambda\,c_{{0}} \left( 36\,{c_{{0}}}^{2}b_{{0}}+341\,\lambda-72\,{b_{{0}}}^{2}d_{{0}}-36\,\delta_{{3}}{b_{{0}}}^{2}
 \right) {a_{{1}}}^{2}}{ 1152 \, {b_{{0}}}^{8}}}\\
&&-{\frac { \left( \lambda\,c_{{0}}c_{{1}}+\lambda\,d_{{1}}b_{{0}}+60\,g_{{1}}{b_{{0}}}^{4}+48\,\gamma_{{1}}{b_{{0}}}^{4} \right) a_{{1}}}{8 \, {b_{{0}}}^{5}}}\\
&&+{\frac {\lambda\, \left( -72\,{b_{{0}}}^{2}d_{{0}}+25\,\lambda+36\,{c_{{0}}}^{2}b_{{0}}-36\,\delta_{{3}}{b_{{0}}}^{2} \right) a_{{2}}}{ 288 \, {b_{{0}}}^{6}}}+{\frac {\lambda\,c_{{2}}}{12 \, {b_{{0}}}^{3}}}.
\end{eqnarray*}

\end{appendix}

\end{document}